\documentclass{amsart}
\usepackage[utf8]{inputenc}
\usepackage{amsmath,amsfonts,amssymb}
\usepackage{tikz}
\usepackage{float}
\usetikzlibrary{arrows.meta}
\usepackage{datetime}
\usepackage{bm}
\usepackage{url}
\usepackage{hyperref}

\newtheorem{theorem}{Theorem}[section]
\newtheorem{lemma}[theorem]{Lemma}
\newtheorem{proposition}[theorem]{Proposition}

\theoremstyle{definition}
\newtheorem{definition}[theorem]{Definition}

\theoremstyle{remark}
\newtheorem{remark}[theorem]{Remark}

\numberwithin{equation}{section}

\newcommand{\N}{\mathbb{N}}
\newcommand{\E}{\mathbb{E}}

\newcommand{\R}{\mathbb{R}}
\newcommand{\C}{\mathbb{C}}

\newcommand{\abs}[1]{\left|#1\right|}

\newcommand{\ZZ}{\mathbb{Z}}

\newcommand{\half}{\frac{1}{2}}

\def\Res{\mathsf{Res}}
\def\Ind{\mathsf{Ind}}

\usepackage[utf8]{inputenc}
\usepackage{amsmath,amsfonts,amssymb}
\usepackage{tikz}
\usepackage{datetime}
\usepackage{bm}
\usepackage{url}
\usepackage{float}
\usepackage{amsmath,amssymb,amsthm,amscd,eucal}

\usepackage{tikz}
\usetikzlibrary{backgrounds}
\usepackage{genyoungtabtikz}
\font \nmm = cmbx9

\def\orbita|bc|d{{
\begin{tikzpicture}[xscale=.5,yscale=.5,line width=1.25pt] 
\foreach \i in {1,2} 
{ \path (\i,1.25) coordinate (T\i); \path (\i,.25) coordinate (B\i); } 
\filldraw[fill= black!12,draw=black!12,line width=4pt]  (T1) -- (T2) -- (B2) -- (B1) -- (T1);
\draw[blue] (T1) -- (B2); 
\foreach \i in {1,2} 
{ \filldraw[fill=white,draw=black,line width = 1pt] (T\i) circle (4pt); \filldraw[fill=white,draw=black,line width = 1pt]  (B\i) circle (4pt); } 
\end{tikzpicture}}}

\theoremstyle{definition}

\theoremstyle{remark}

\numberwithin{equation}{section}

\newlength\cellsize 

\savebox2{
\begin{picture}(16,16)
\put(0,0){\line(1,0){16}}
\put(0,0){\line(0,1){16}}
\put(16,0){\line(0,1){16}}
\put(0,16){\line(1,0){16}}
\end{picture}}

\newcommand\boxify[1]{\def\thearg{#1}\def\nothing{}%
\ifx\thearg\nothing\vrule width0pt height\cellsize depth0pt%
  \else\hbox to 0pt{\usebox2\hss}\fi%
  \vbox to \cellsize{\vss\hbox to \cellsize{\hss$_{#1}$\hss}\vss}}

\savebox3{
\begin{picture}(8,8)
\put(4,4){\circle{8}}
\end{picture}}

\newcommand{\circify}[1]{\def\thearg{#1}\def\nothing{}%
\ifx\thearg\nothing\vrule width0pt height\cellsize depth0pt%
  \else\hbox to 0pt{\usebox3\hss}\fi%
  \vbox to \cellsize{\vss\hbox to \cellsize{\hss$_{#1}$\hss}\vss}}

\newcommand\nullify[1]{\def\thearg{#1}\def\nothing{}%
\ifx\thearg\nothing\vrule width0pt height\cellsize depth0pt%
  \else\hbox to 0pt{\hss}\fi%
  \vbox to \cellsize{\vss\hbox to \cellsize{\hss$_{#1}$\hss}\vss}}

\newcommand\tableau[1]{\vtop{\let\\=\cr
\setlength\baselineskip{-16000pt}
\setlength\lineskiplimit{16000pt}
\setlength\lineskip{0pt}
\halign{&\boxify{##}\cr#1\crcr}}}

\newcommand\cirtab[1]{\vline\vtop{\let\\=\cr
\setlength\baselineskip{-8000pt}
\setlength\lineskiplimit{8000pt}
\setlength\lineskip{0pt}
\halign{&\circify{##}\cr#1\crcr}}}

\newcommand\nulltab[1]{\vtop{\let\\=\cr
\setlength\baselineskip{-8000pt}
\setlength\lineskiplimit{8000pt}
\setlength\lineskip{0pt}
\halign{&\nullify{##}\cr#1\crcr}}}

\title[Cycle statistics and representation theory]{Representation theory and cycle statistics for random walks on the symmetric group}
\author{Dominic Arcona}
\address{Department of Mathematics, University of Southern California, Los Angeles, CA 90089-2532, USA }
\email{arcona@usc.edu}
\date{August 27, 2026.} 

\begin{document}

\begin{abstract}
We use representation theory of $S_n$ to analyze the mixing of cycle type statistics $a_j(\sigma) = \{$\# of $j$-cycles of $\sigma$\} for any fixed $j$ in permutations $\sigma_t$ resulting from the $t$-step random transposition walk on $S_n$. We also derive analogous results for the star transposition walk. Our approach uses the method of moments; a key ingredient is a new formula for the coefficients in the irreducible character decomposition of the $S_n$-class function $(a_j)^r(\sigma)=\{(\text{\# of $j$-cycles of $\sigma$})^r\}$ for any positive integers $r,j$ when $n\geq 2rj$.
\end{abstract}
\maketitle
\tableofcontents
\section{Introduction}
In the symmetric group $S_n$, it is well known that for any fixed positive integer $j$, the number of $j$-cycles of a uniformly random permutation converges to a Poisson$(1/j)$ random variable as $n\rightarrow \infty$; see for instance \cite{AT92}. In this paper, we examine the distribution of $j$-cycles of random permutations obtained by sampling from the $t$-fold convolution of a probability measure defining one of two simple random walks: the \textit{star transposition walk} or the \textit{random transposition walk}, both defined in Section \ref{prelim}. We leverage the ``method of moments'' \cite{D87}, so it is useful to know that a Poisson($\lambda$) random variable has $r$th moment equal to \[
\sum_{a=0}^r S(r,a) \lambda^a, 
\]
where $S(r,a)$ is the Stirling number of the second kind, the number of set partitions of $[r]$ into $a$ blocks. 

Let us clarify some notation on partitions before proceeding. We say that a weakly decreasing sequence of nonnegative integers $\lambda = (\lambda_1, \lambda_2, \dots)$ denotes a partition of size $|\lambda|$. We write $\lambda_i$ for the size of the $i$th part of $\lambda$, and $d_{\lambda}$ for the number of standard Young tableaux of shape $\lambda$. We write $\lambda \vdash n$ when $\abs{\lambda} = n$ and frequently identify $\lambda$ with its Young diagram. If $\lambda = (\lambda_1,\lambda_2,\ldots)$, we use $\lambda^\prime$ to denote the partition obtained by transposing rows and columns of $\lambda$ and $\bar{\lambda}$ to refer to the partition $(\lambda_2,\lambda_3,\ldots)$ of size $|\bar{\lambda}| = |\lambda| - \lambda_1$ so that $\lambda = (\lambda_1, \bar{\lambda})$. We say $\sigma \in S_n$ has cycle type $\alpha = (1^{m_1}2^{m_2}\cdots n^{m_n})$ if $\sigma$ has $m_j$ $j$-cycles for each $j$ between 1 and $n$. The cycle type of a permutation determines which conjugacy class of $S_n$ it belongs to. We will write partitions $\alpha \vdash n$ indexing conjugacy classes in decreasing order of cycle length; for example, the conjugacy class of a transposition is indexed by the partition $(2,1^{n-2})$.   

Our approach will use representation theory of $S_n$ to compute the moments of the class function 
 \[
 a_j(\sigma) = \{\text{\# of $j$-cycles of $\sigma$}\}
 \]
 over probability measures of interest. This analysis uses Theorem \ref{thm:A}, our central algebraic result, which gives the irreducible character decomposition of the class functions 
 \[
 (a_j)^r(\sigma)=\{(\text{\# of $j$-cycles of $\sigma$})^r\}
 \]
 for any positive integers $r,j$ when $n\geq 2rj$. 

\begin{theorem}\label{thm:A}
    Let $r,n \in \N$, $j \in \ZZ_{\geq 1}$, and $\sigma \in S_n$ with $n\ge 2rj$. Let 
    \[
    (a_j)^r(\sigma) = \{(\text{\#$j$-cycles of $\sigma$})^r\}.
    \] Then we have the irreducible $S_n$-character decomposition:
    \[(a_j)^r = \frac{1}{j^r}\left(\sum_{\lambda \vdash n} m_{\lambda,r}^j \chi^\lambda\right).\]
    The (signed) integer multiplicity $m_{\lambda,r}^j$ is 0 unless $\lambda=(n-tj,\bar{\lambda})$ for $t\in \N_{\leq r}$. For $\lambda=(n-tj,\bar{\lambda})$, it is given by
    \[
    m^{j}_{\lambda,r}
    =
    \chi^{\bar{\lambda}}_{(j^t)}\sum_{a=t}^r S(r,a)\binom{a}{t}j^{r-a},
    \]
    where $\chi^{\bar{\lambda}}_{(j^t)}$ is the irreducible $S_{tj}$-character indexed by $\bar{\lambda}$ evaluated on conjugacy class $(j^t)$ and $S(r,a)$ denotes the Stirling number of the second kind.
\end{theorem}

This decomposition is new for $r>1$ except in the case that $j=1$ (fixed points), where it is known by several methods \cite{BHH16,D16,F24}. Its proof is inspired by the restriction-induction Bratteli diagrams of Benkart, Halverson, and Harmon \cite{BHH16}, discussed in Section 6. The irreducible character decomposition of $a_j$ is given by Alon and Kozma in \cite{AK13} and stated in Section \ref{reptheory}. Additional context and representation-theoretic background for Theorem \ref{thm:A} is also given in Section \ref{prelim}.

We obtain three main results concerning the limiting distribution of cycle type statistics in the star transposition and random transposition walk settings by applying Theorem \ref{thm:A}. Both walks are defined in Section \ref{rwalks}, where an overview of a rich literature of related results can also be found. We believe that Theorem \ref{thm:A} can be applied to analyze finite $j$-cycle mixing in other random walks on $S_n$ such as the random $i$-cycle walk as well.

The first application establishes a Poisson limiting distribution for the number of fixed points in permutations obtained from $t = \lfloor n\log n + cn\rfloor$ applications of star transpositions. Theorem \ref{thm:fixstar} requires only the known $j=1$ specialization of Theorem \ref{thm:A} given previously by each of \cite{BHH16,D16,F24}.
\begin{theorem}\label{thm:fixstar}
    For $c$ a fixed real number, the distribution of the number of fixed points after multiplying $\lfloor n\log n + cn\rfloor$ star transpositions starting from the identity converges to a Poisson$(1+e^{-c})$ limit as $n\rightarrow \infty$.
\end{theorem}
This result parallels a recent result of Fulman \cite{F24}, who used the method of moments to prove a similar limiting Poisson distribution for fixed points of permutations obtained via $i$-cycle shuffling. Since the total variation limit profile for this walk established by Nestoridi \cite{N22} coincides with the same timescale, $\lfloor n \log n + cn\rfloor$, Theorem \ref{thm:fixstar} confirms that fixed points are among the last permutation statistics to mix for this walk.

The second application, Theorem \ref{thm:B}, shows that, for $j\geq 2$, the finite $j$-cycles mix considerably faster than fixed points in the star transposition setting, the critical timescale being just $\lfloor cn\rfloor$ many steps with $c$ a positive real number. 
\begin{theorem}\label{thm:B}
     Fix an integer $j\geq 2$ and a real number $c>0$. The distribution of the number of $j$-cycles after multiplying $\lfloor cn\rfloor$ star transpositions starting from the identity converges to a Poisson$\left(\frac{1}{j}(1-e^{-jc})\right)$ limit as $n\rightarrow \infty$. 
\end{theorem}

Our final result describes the limiting distribution of finite $j$-cycles and their critical timescale for mixing in the random transposition walk setting. 
\begin{theorem}\label{thm:rtransp}
    Fix a positive integer $j$ and a real number $c$. The distribution of the number of $j$-cycles after multiplying $\lfloor \frac{n}{2j}(\log n + (j-1)\log\log n + c)\rfloor$ random transpositions starting from the identity converges to a Poisson$\left(\frac{1}{j}(1+\frac{e^{-c}}{j!})\right)$ limit as $n\rightarrow \infty$. 
\end{theorem}
Theorem \ref{thm:rtransp} generalizes a result of Matthews \cite{M88} that the fixed point count converges to Poisson($1+e^{-c})$ after $\lfloor \frac{n}{2}(\log n + c)\rfloor$ steps. In addition, it establishes a hierarchy of critical windows $t_1 > t_2 > \cdots > t_j = \lfloor \frac{n}{2j}(\log n + (j-1)\log\log n + c)\rfloor$ for the convergence of finite-cycle statistics to their stationary Poisson laws.

Let us now give some probabilistic intuition for Theorem \ref{thm:rtransp}. We can see that at least the mean number of $j$-cycles claimed at time $t_j$ is approximately correct using the Erd\H{o}s--R\'enyi random graph associated with the random transposition walk. Concretely, let $G_t = (V, E_t)$ with vertex set $V = [n]$ and edge set $E_t$ containing edge $(i,j)$ if transposition $(i\,j)$ has been selected and applied in the first $t$ steps of the random transposition walk. Then $G_{t} \approx G(n,p_{t})$ with 
\[
p_t = 1 - \left(1-\binom{n}{2}^{-1}\right)^t,
\]
where $G(n,p)$ denotes the Erd\H{o}s--R\'enyi random graph on $n$ vertices with probability $p$ of any given edge being present. The Erd\H{o}s--R\'enyi theorem then says that after $t \sim cn$ steps with $c>1/2$, a macroscopic connected component $\mathcal{C}_t$ of size $\Theta(n)$ exists in $G_t$ with high probability. A conjecture of Berestycki \cite[Conjecture 1.3]{T20} asserts that permutations obtained from $t$ steps of the random transposition walk restricted to positions labeled by vertices in $\mathcal{C}_t$ are approximately uniformly distributed in $S_{\abs{\mathcal{C}_t}}$. Thus at time $t_j \sim \frac{n}{2j}\log n \gg \frac{n}{2}$, we expect $1/j$ finite $j$-cycles to come from $\mathcal{C}_{t_j}$. 

Outside of $\mathcal{C}_{t_j}$, we reason as follows. $G_{t_j} \approx G(n,p_{t_j})$ with $p_{t_j} = 1-\mathrm{exp}(-t_j/\binom{n}{2})$. Because $p_{t_j}\rightarrow 0$ as $n$ grows, the probability that any connected component of size $j$ in $G(n,p_{t_j})$ is a tree of size $j$ is asymptotically 1. Now let $\E_{t_j}[T_j]$ denote the expected number of trees of size $j$ in $G_{t_j}$. One can check that, asymptotically, 
\[
\E_{t_j}[T_j] \sim \frac{e^{-c}}{j \cdot j!},
\]
which is already known from the theory of random graphs \cite[Theorem 5.4]{B01}. Since repeated transpositions occurring within connected components of size $j$ are very unlikely (by the same reasoning that finite-sized connected components are almost all trees), these trees of size $j$ each contribute $(1-o_n(1))$ $j$-cycles, bringing the total expected number of $j$-cycles in a permutation obtained from a random transposition walk of length $t_j$ steps to $\frac{1}{j}+\frac{e^{-c}}{j\cdot j!}-o(1)$.   

The remainder of this paper is organized as follows. Section \ref{prelim} provides useful facts from the representation theory of $S_n$ and symmetric function theory in \ref{reptheory} and \ref{subsec: symfuns}; it concludes with background on the random transposition and star transposition walks including a framework for their representation-theoretic analysis in \ref{rwalks}. Theorem \ref{thm:fixstar} is proved in Section \ref{starwalk}, which examines the mixing of fixed points after repeated application of star transpositions. Section \ref{limj} treats finite $j$-cycle mixing in the star transposition setting and proves Theorem \ref{thm:B}. Section \ref{rtransp} is dedicated to the mixing of finite $j$-cycles in the random transposition walk setting and proves Theorem \ref{thm:rtransp}. Finally, Section \ref{bratteli} contains proof of Theorem \ref{thm:A} using combinatorial tools.
 
\section{Preliminaries}\label{prelim}

\subsection{Representation theory of \(S_n\)}\label{reptheory}\hfill\\

In this section, we introduce some basic facts from the representation theory of $S_n$ that will be useful for our analysis. For a more thorough reference, see Sagan \cite{S01} or James and Kerber \cite{JK81}. 

The irreducible representations of $S_n$ over $\mathbb{C}$ (isomorphic to Specht modules) are indexed by partitions of \( n \), and a classical result (Maschke's Theorem) says that every finite-dimensional representation of \( S_n \) over \( \mathbb{C} \) decomposes into a direct sum of irreducible representations. We let $S^\lambda$ denote the irreducible representation of $S_n$ indexed by $\lambda\vdash n$ and $\chi^\lambda$ denote its character, i.e. $\chi^\lambda(\sigma) = \mathrm{Tr}(S^\lambda(\sigma))$.  With $\lambda,\, \alpha \vdash n$, we write $\chi^\lambda_\alpha$ to denote the value the irreducible $S_n$-character indexed by $\lambda$ evaluated on conjugacy class $\alpha$. We assume familiarity with the Murnaghan-Nakayama rule, stated below.
\begin{theorem}\cite[Theorem 4.10.2]{S01}\label{MN}
    If $\lambda \vdash n$ and $\alpha = (\alpha_1,\ldots,\alpha_k)$ is a composition of $n$, then 
    \[
    \chi^\lambda_\alpha = \sum_{\substack{\mu\vdash n-\alpha_1 \\ \lambda/\mu \textrm{ a rim-$\alpha_1$ hook}}} (-1)^{ll(\lambda/\mu)}\chi^{\mu}_{\bar{\alpha}},
    \]
    where $\bar{\alpha} = \alpha/\alpha_1 = (\alpha_2,\ldots,\alpha_k)$ and a rim-$\alpha_1$ hook is defined as a rim-hook having $\alpha_1$ cells.
\end{theorem}

Identifying $\lambda$ with its Young diagram, the dimension of the irreducible representation \( S^\lambda \) is equal to $d_\lambda$ (equivalently $\chi^\lambda_{(1^n)}$), given by the famous hook-length formula of Frame-Robinson-Thrall \cite{FRT54}:
\begin{equation}\label{eqn:hlf}
    \dim S^\lambda = d_\lambda = \frac{n!}{\prod_{\square \in \lambda} h(\square)}
\end{equation}

where \( h(\square) \) is the hook-length of a box \( \square \) in \( \lambda \).  

A function $f: G \rightarrow \C$ is called a \textit{class function} if it is constant on conjugacy classes of $G$. The irreducible characters of $S_n$ form a basis for the complex vector space
$\mathrm{CF}(S_n)$ of class functions on $S_n$.  Let $R(S_n)$ denote the Grothendieck
group of finite-dimensional $S_n$-representations, so that
$R(S_n)\otimes_{\mathbb{Z}}\mathbb{C}$ consists of complex linear combinations
of irreducible representations (irreps for short).  The character map induces an isomorphism
\[
R(S_n)\otimes_{\mathbb{Z}}\mathbb{C} \;\xrightarrow{\;\sim\;}\; \mathrm{CF}(S_n),
\]
and therefore every class function can be viewed as the character of a (complex) virtual
$S_n$-representation.
 We will consider class functions such as $a_j(\sigma) = $ \{\# of $j$-cycles of $\sigma$\} and the corresponding virtual representations $\rho_{a_j}$. 

As an example, let $\varrho_n$ denote the $n$-dimensional defining representation of $S_n$, which decomposes into irreps as
\[\varrho_n \cong S^{(n)} \oplus S^{(n-1,1)}.\]
At any $\sigma\in S_n$ the trace of $\varrho_n(\sigma)$ is exactly the number of fixed points of $\sigma$, which makes it identical to $\rho_{a_1}$, which has the class function $a_1(\sigma)$ = \{\# of fixed points of $\sigma$\} as its character. Thus we have the irreducible character decomposition: \[a_1 = \chi^{(n)} + \chi^{(n-1,1)}.\]

Alon and Kozma \cite{AK13} provide an irreducible character decomposition of $a_j$ in general; a specialization to the case $n \geq 2j$ is given as Theorem \ref{thm:10} below.

\begin{theorem}\cite[Theorem 3]{AK13}\label{thm:10}
    For $1\leq j \leq \frac{n}{2}$,
    \begin{equation}\label{eqn:1}
    a_j = \frac{1}{j}\left(1+\sum_{i=0}^{j-1}(-1)^i\chi^{(n-j,j-i,1^i)}\right).
    \end{equation}
\end{theorem}

Now let $\psi_j(\sigma)$ denote the class function $j\,a_j(\sigma)=$ \{$j$\,(\#$j$-cycles of $\sigma$)\} and let $\rho_{\psi_j}$ be the (unique) virtual $S_n$-representation with character $\psi_j$. Translating Equation ~\eqref{eqn:1} from characters to representations and letting `+' and `$-$' denote formal sum and difference gives, for $n\geq 2j$,
\begin{equation}
    \rho_{\psi_j} \cong S^{(n)} + \sum_{i=0}^{j-1}(-1)^i{S^{(n-j,j-i,1^i)}}.
\end{equation}

One reason we consider $\rho_{\psi_j}$ instead of $\rho_{a_j}$ is that the signed multiplicities of irreps in the irreducible decomposition of $\rho_{\psi_j}^{\otimes r}$ are integers. These are precisely the coefficients $m_{\lambda,r}^j$ appearing in Theorem \ref{thm:A}. We will see in Section \ref{bratteli} how to compute tensor powers of $\rho_{\psi_j}$ using a restriction-induction Bratteli diagram, which in turn affords us the ability to express the class function $(a_j)^r(\sigma) = $ \{(\#$j$-cycles of $\sigma)^r$\} as a linear combination of irreducible characters for any $r,j$ when $n$ is large enough.

\subsection{Symmetric and shifted-symmetric functions}\label{subsec: symfuns}\hfill\\

We briefly recall some aspects of symmetric function theory used in Sections \ref{rtransp} and \ref{bratteli}. Standard references are Macdonald~\cite[Chapter 1]{M95} and Sagan~\cite[Chapter 4]{S01}. For additional background on shifted-symmetric functions, see Okounkov and Olshanski~\cite{OO97} or Okounkov and Pandharipande~\cite{OP06}.

Let
\[
\Lambda=\bigoplus_{n\ge0}\Lambda_n
\]
denote the graded ring of symmetric functions over \(\mathbb{Q}\), where \(\Lambda_n\) is the space of homogeneous symmetric functions of degree \(n\). For a partition
\[
\alpha=(1^{m_1}2^{m_2}\cdots)\vdash n,
\]
define
\[
z_\alpha=\prod_{i\ge1}i^{m_i(\alpha)}m_i(\alpha)!.
\]
The power-sum symmetric functions are
\[
p_k=\sum_{i\ge1}x_i^k,
\qquad
p_\alpha=\prod_i p_{\alpha_i},
\]
and \(\{p_\alpha:\alpha\vdash n\}\) forms a basis of \(\Lambda_n\). Another distinguished basis is the Schur basis
\[
\{s_\lambda:\lambda\vdash n\},
\]
which is orthonormal with respect to the Hall inner product,
\[
\langle s_\lambda,s_\mu\rangle=\delta_{\lambda\mu},
\]
or equivalently,
\[
\langle p_\alpha,p_\beta\rangle
=
z_\alpha\,\delta_{\alpha\beta}.
\]

The transition matrix between the Schur and power-sum bases is given by the irreducible characters of the symmetric group. Writing \(\chi^\lambda_\alpha\) for the irreducible character of \(S_n\) indexed by \(\lambda\vdash n\), evaluated on the conjugacy class of cycle type \(\alpha\vdash n\),
\begin{equation}\label{eq:schur-powersum}
s_\lambda
=
\sum_{\alpha\vdash n}
\frac{\chi^\lambda_\alpha}{z_\alpha}p_\alpha,
\qquad
p_\alpha
=
\sum_{\lambda\vdash n}
\chi^\lambda_\alpha s_\lambda.
\end{equation}
This correspondence is encoded by the \textit{Frobenius characteristic map}
\[
\operatorname{ch}:\mathrm{CF}(S_n)\longrightarrow\Lambda_n,
\qquad
\operatorname{ch}(f)
=
\sum_{\alpha\vdash n}
\frac{f(\alpha)}{z_\alpha}p_\alpha,
\]
where \(\mathrm{CF}(S_n)\) denotes the space of class functions on \(S_n\). In particular,
\[
\operatorname{ch}(\chi^\lambda)=s_\lambda.
\]
Moreover, \(\operatorname{ch}\) is an isometry between the standard inner product on class functions,
\[
\langle f,g\rangle_{S_n}
=
\frac1{n!}\sum_{\sigma\in S_n}f(\sigma)g(\sigma)
=
\sum_{\alpha\vdash n}\frac{f(\alpha)g(\alpha)}{z_\alpha},
\]
and the Hall inner product on \(\Lambda_n\). Consequently, the orthogonality relations
\[
\sum_{\alpha\vdash n}
\frac{\chi^\lambda_\alpha\chi^\mu_\alpha}{z_\alpha}
=
\delta_{\lambda\mu},
\qquad
\sum_{\lambda\vdash n}
\chi^\lambda_\alpha\chi^\lambda_\beta
=
z_\alpha\delta_{\alpha\beta},
\]
translate precisely into the orthogonality of the Schur and power-sum bases.

Many operators arising in algebraic combinatorics are diagonal in the Schur basis. If
\[
T:\Lambda_n\rightarrow\Lambda_n,
\qquad
T(s_\lambda)=\theta_\lambda s_\lambda,
\]
then the eigenvalue function
\[
\lambda\longmapsto\theta_\lambda
\]
completely determines \(T\). Throughout this paper we freely pass between Schur-diagonal operators and their eigenvalue functions on partitions.

Unlike ordinary symmetric functions, these eigenvalue functions naturally depend on the shifted row coordinates of a Young diagram. To study such eigenvalue functions, it is natural to work in the algebra of shifted symmetric functions, denoted $\Lambda^*$. Following Okounkov and Olshanski~\cite{OO97}, \(\Lambda^*\) is the filtered algebra obtained as the inverse limit of shifted symmetric polynomials, namely polynomial functions in the shifted row coordinates
\[
x_i=\lambda_i-i
\]
that are symmetric in the variables \(x_1,x_2,\ldots\) and satisfy the usual stability condition. Unlike the ordinary symmetric-function ring \(\Lambda\), the algebra \(\Lambda^*\) is naturally \emph{filtered}, not graded. The filtration is inherited from polynomial degree in the shifted coordinates:
\[
\Lambda^*_{\le0}
\subseteq
\Lambda^*_{\le1}
\subseteq
\Lambda^*_{\le2}
\subseteq
\cdots,
\]
where \(\Lambda^*_{\le d}\) consists of shifted symmetric functions of shifted degree at most \(d\).

A particularly useful family of shifted symmetric functions is given by the shifted complete homogeneous symmetric functions \(h_r^*\), defined by the generating series
\begin{equation}\label{eq:shifted-generating}
\prod_{i\ge1}\frac{u+i}{u+i-\lambda_i}
=
\sum_{r\ge0}\frac{h_r^*(\lambda)}{u^{\downharpoonright r}},
\qquad
u^{\downharpoonright r}
=
u(u-1)\cdots(u-r+1),
\end{equation}
which appears as Equation~(12.3) of~\cite{OO97}. Okounkov and Olshanski prove that \(h_r^*\) has shifted degree \(r\), that is,
\[
h_r^*\in\Lambda^*_{\le r},
\]
and that the family \(\{h_r^*\}_{r\ge0}\) generates the filtered algebra \(\Lambda^*\). The generating series \eqref{eq:shifted-generating} will play a central role in the proof of Lemma~\ref{lemma:dim}.

Another important generating family is given by the shifted power sums
\begin{equation}\label{eqn:powersum}
p_k^*(\lambda)
=
\sum_{i\ge1}
\left[
\left(\lambda_i-i+\frac12\right)^k
-
\left(-i+\frac12\right)^k
\right],
\qquad
k\ge1,
\end{equation}
each of which has shifted degree \(k\). 

The main role of \(\Lambda^*\) in this paper is to organize eigenvalue functions of Schur-diagonal operators. Whenever an eigenvalue function can be recognized as a polynomial function of the shifted row data of a Young diagram $\lambda$, it determines an element of \(\Lambda^*\).

An example meriting further discussion is the \textbf{cut-and-join operator} on $\Lambda$, defined by 
\begin{equation}\label{def:cut-and-join}
 \Delta = \frac12\sum_{a,b\ge1}(a+b)p_ap_b\,\frac{\partial}{\partial p_{a+b}} +
\frac12\sum_{a,b\ge1}ab\,p_{a+b}\,\frac{\partial}{\partial p_a}\frac{\partial}{\partial p_b}.   
\end{equation}

Since \(\Delta\) preserves homogeneous degree, it restricts to an endomorphism of each \(\Lambda_n\). Under the Frobenius characteristic map, left multiplication by the transposition class sum
\[
C_{(2)}
=
\sum_{1\le i<j\le n}(ij)
\]
in the center \(Z(\mathbb C[S_n])\) corresponds precisely to the operator \(\Delta\). Since \(C_{(2)}\) is central, it acts by a scalar on every irreducible representation \(S^\lambda\), and consequently \(\Delta\) acts diagonally on the Schur basis:
\[
\Delta s_\lambda
=
\frac{\kappa(\lambda)}2\,s_\lambda,
\]
where
\[
\kappa(\lambda)
=
\sum_i
\lambda_i(\lambda_i-2i+1)
\]
is a shifted symmetric function of shifted degree \(2\).

The terminology ``cut-and-join'' reflects the combinatorics of multiplying a permutation by a transposition. If the transposition exchanges two elements lying in different cycles, those cycles merge into a single cycle; if it exchanges two elements lying in the same cycle, that cycle splits into two. Correspondingly, writing
\[
\Delta=\Delta_++\Delta_-,
\]
where
\[
\Delta_+
=
\frac12
\sum_{a,b\ge1}
(a+b)p_ap_b\,\frac{\partial}{\partial p_{a+b}},
\qquad
\Delta_-
=
\frac12
\sum_{a,b\ge1}
ab\,p_{a+b}\,\frac{\partial}{\partial p_{a}}\frac{\partial}{\partial p_{b}},
\]
the operator \(\Delta_+\) replaces one factor $p_{a+b}$ 
by the product $p_a p_b$, thereby increasing the number of parts of the indexing partition by one, while \(\Delta_-\) similarly decreases the number of parts by one. Thus the two differential operators exactly mirror the two possible outcomes of multiplying a permutation by a transposition.

This combinatorial interpretation will play an essential role in Section~\ref{rtransp}. Although the eigenvalue function \(\kappa(\lambda)\) has shifted degree \(2\), multiplication by the transposition class sum changes the number of cycles in a conjugacy class by at most one. Consequently, if
\[
\mathsf P_{\le\ell}
=
\operatorname{span}\{p_\nu:\ell(\nu)\le\ell\},
\]
then 
\[
\Delta(\mathsf P_{\le\ell})
\subseteq
\mathsf P_{\le\ell+1}.
\]

This fundamental property of the cut-and-join operator is used throughout Section~\ref{rtransp}.

\subsection{Evolution of cycle type statistics over random walks}\label{rwalks}\hfill\\

In this subsection, we make clear how one can make use of representation theoretic tools to analyze the distribution of a statistic of interest as it evolves over the course of a random walk on $S_n$.

\begin{definition}\label{def:ftrans}
    Let $P$ be a probability on a group $G$. The Fourier transform of $P$ at representation $\rho$ is defined as $\widehat{P}(\rho) = \sum_{g\in G}P(g)\rho(g)$.
\end{definition}

Suppose we have a random walk on $G$ governed by probability law $P$. Concretely, we mean that at each step we begin at element $g \in G$, sample an element $h \in G$ according to $P$, and proceed to element $hg \in G$. When we take multiple steps, say $t$ of them, each according to $P$, the result is equivalent to sampling $h^\prime \in G$ according to the $t$-fold convolution measure $P^{*t}$ and stepping from $g$ to $h^\prime g$. A convenient property of the Fourier transform (see Exercise 1 in \cite{Dia88}) is that
\[\widehat{P^{*t}}(\rho) = \widehat{P}(\rho)^t.\]  

A natural way to study the evolution of a statistic over a random walk is to examine the moments of its distribution after $t$ steps. Theorem \ref{thm:1}, used implicitly by Fulman \cite{F24}, Diaconis \cite{D88}, and others, explains a key connection that makes representation theory useful in the study of random walks on $S_n$.

\begin{theorem}\label{thm:1}
    Let $t\in \N$, $P$ a probability on $S_n$ with Fourier transform $\widehat{P}$, and $\upsilon$ an $S_n$-class function with irreducible character decomposition
\[\upsilon = \sum_{\lambda\vdash n}c_\lambda \chi^\lambda,\]
where the $c_\lambda$ are $\mathbb{C}$-valued constants.
The expectation of $\upsilon$ after $t$ steps of a random walk on $S_n$ governed by $P$ starting from the identity permutation is given by
\[\E_{P^{*t}}[\upsilon] = \sum_{\lambda\vdash n}c_\lambda\mathrm{Tr}\left(\widehat{P}(S^{\lambda})^t\right).\]
\end{theorem}

\begin{proof}
\begin{alignat*}{2}
\E_{P^{*t}}[\upsilon] &= \sum_{\pi\in S_n}P^{*t}(\pi)\upsilon(\pi)\\
    &= \sum_{\pi}P^{*t}(\pi)\sum_{\lambda \vdash n}c_\lambda\chi^{\lambda}(\pi)\\
    &= \sum_{\lambda}c_\lambda\sum_{\pi}P^{*t}(\pi)\mathrm{Tr}(S^{\lambda}(\pi))\\
    &= \sum_{\lambda}c_\lambda\mathrm{Tr}\left(\sum_{\pi}P^{*t}(\pi)S^{\lambda}(\pi)\right)\\
    &= \sum_{\lambda}c_\lambda\mathrm{Tr}\left(\widehat{P^{*t}}(S^{\lambda})\right)\\
    &= \sum_{\lambda}c_\lambda\mathrm{Tr}\left(\widehat{P}(S^{\lambda})^t\right). 
\end{alignat*}
\end{proof}

This gives us a highly flexible tool for analyzing random walks on $S_n$; for any statistic, shuffling law pair $(\upsilon,P)$, we can compute $\E_{P^{*t}}[\upsilon]$ for any $t$ given only the irreducible character decomposition of $\upsilon$ and some knowledge of $\widehat{P}(S^\lambda)$. The irreducible character decomposition of class function $(a_j)^r(\sigma) = \{(\text{\# of $j$-cycles of $\sigma$})^r\}$ will be a main concern of Section \ref{bratteli}. Let us conclude this section with a discussion of the two shuffles we analyze, including expressions for their respective probability law Fourier transforms at any given irrep $S^\lambda$. Note that we use \textit{walk} and \textit{shuffle} interchangeably to describe permutation mixing processes throughout the paper.  

\begin{definition}\label{def:star}
The \textbf{star transposition walk} is defined by the probability measure $Q$ on $S_n$ as follows.
    \[Q(\pi) = 
    \begin{cases} 
      \frac{1}{n-1} & \text{if $\pi$ is one of the transpositions $(1i) \qquad 2\leq i \leq n$.}\\
      0 & \text{otherwise} 
   \end{cases}
\]
\end{definition}

The star transposition shuffle, a natural variant of the random transposition shuffle, was first analyzed by Flatto, Odlyzko, and Wales \cite{FOW85}, who found the eigenvalues of the associated Markov chain using representation theory. Diaconis \cite{D88} subsequently proved a mixing time of $\lfloor n\log n +cn\rfloor$ steps, also using representation theory. He also found the mean and variance of the distribution of fixed points after applying $\lfloor n\log n + cn\rfloor$ star transpositions both converge to $1+e^{-c}$ for fixed $c$ as $n \rightarrow \infty$. More recently, Nestoridi \cite{N22} proved the total variation limit profile.

We follow the convention of \cite{FOW85}, excluding the identity step as a valid star transposition, while \cite{D88} and \cite{N22} allow it. The eigenvalues of the Markov chain associated with the random walk are nearly equivalent in either case and the asymptotics of the walk are unaffected by the exclusion.

It is notable that $Q$ is non-constant on conjugacy classes, yet representation theory is still useful in our analysis. The reason is that $\widehat{Q}(S^\lambda)$ is diagonal in a suitable basis \cite{D88}, and we can still find a compact expression for $\mathrm{Tr}(\widehat{Q}(S^\lambda)^t)$. The following Lemma, adapted from \cite{D88}, gives such an expression.

\begin{lemma}\label{lemma:201}
Let $Q$ be the probability law on $S_n$ given by definition \ref{def:star} with Fourier transform $\widehat{Q}$, $S^\lambda$ the irrep indexed by $\lambda \vdash n$, and $t \in \N$. Then
\[\mathrm{Tr}\left(\widehat{Q}(S^\lambda)^t\right) = \sum_{\lambda^i}d_{\lambda^i}\left(\frac{\lambda_i-i}{n-1}\right)^t,\]

\noindent where $\lambda^i$ denotes the partition of $n-1$ obtained by removing an inner corner from the $i$th row of $\lambda$, and $\lambda_i$ denotes the length of the $ith$ row of $\lambda$. The sum ranges over all valid inner corner removals from $\lambda$, equivalently over all $i$ such that $\lambda^i$ has a valid Young diagram.
\end{lemma}

We now turn our attention to the random transposition walk on $S_n$.

\begin{definition}\label{def:icycle}
The \textbf{random transposition walk} is defined by probability measure $P$ on $S_n$ as follows.
    \[P(\pi) = 
    \begin{cases} 
      \frac{1}{\binom{n}{2}} & \text{if $\pi$ is a 2-cycle}\\
      0 & \text{otherwise} 
   \end{cases}
\]
\end{definition}

A seminal paper of Diaconis and Shashahani \cite{DS81} showed that the walk undergoes a cutoff in total variation mixing after $\frac{n}{2}\log n$ steps with window size $O(n)$. More specifically, they showed that if $t=\lfloor\frac{n}{2}(\log n + c)\rfloor$, then
\[
d_{TV}(P^{*t}, U) \xrightarrow[c \to \infty]{n \to \infty} 0 \quad \text{and} \quad d_{TV}(P^{*t}, U) \xrightarrow[c \to -\infty]{n \to \infty} 1,
\]
where we denote the uniform measure on $S_n$ by $U$. The total variation limit profile,
\[
d_{TV}(P^{*t}, U) = d_{TV}(\mathrm{Poisson}(1+e^{-c}),\mathrm{Poisson}(1)),
\]
was later proved by Teyssier in \cite{T20}. The distribution of the number of fixed points after $t$ steps was shown to be $\mathrm{Poisson}(1+e^{-c})$ by Matthews \cite{M88}. On timescale $t$, the total variation distance to stationarity of the walk mirrors exactly the distance between the fixed point count distribution under $P^{*t}$ and its distribution under the uniform measure $U$, the stationary distribution of the walk. This fact is more than coincidence, as recent results of Jain and Sawhney \cite{JS24} and Teyssier \cite{T26} demonstrate.

One hope is that the analysis of finite $j$-cycle mixing can be extended to other random walks of interest on $S_n$ such as the random $i$-cycle walk, a generalization of the random transposition walk. A Poissonian total variation limit profile \cite{NOT22} and limiting distribution of the fixed point count \cite{F24} were both proved recently in this setting.

The following Lemma will be useful for analyzing $j$-cycle mixing in the random transposition setting. It is well-known; see for instance \cite{Dia88} or \cite{F24}.

\begin{lemma}\label{lemma:202}
Let $P$ be the probability law on $S_n$ given by definition \ref{def:icycle} with Fourier transform $\widehat{P}$, $S^\lambda$ the irrep indexed by $\lambda \vdash n$, and $t \in \N$. Then
\[\mathrm{Tr}\left(\widehat{P}(S^\lambda)^t\right) = d_{\lambda}\left(\frac{\chi^\lambda_\tau}{d_\lambda}\right)^t,\]
where $\chi^\lambda_\tau$ is the character of $S^\lambda$ evaluated on the conjugacy class $(2,1^{n-2})$ of permutations containing a single transposition and $n-2$ fixed points.
\end{lemma}

\section{Fixed points and the star transposition walk}\label{starwalk}

Using a method of moments, Fulman recently proved that the limiting distribution of the number of fixed points after applying $t_i=\lfloor \frac{1}{i}n\log n + cn\rfloor$ random $i$-cycles for $i\geq 2$ and $c \in \R$ fixed converges to a Poisson$(1+e^{-ic})$ distribution as $n\rightarrow \infty$ \cite{F24}. We show that an analogous result holds for the star transposition walk, using $k$ as the step count in place of $t$ in this section to avoid confusion with the number of boxes $t$ below the first row in partitions $\lambda = (n-t,\bar{\lambda})$, which will be germane.

Recall the preliminaries of Section \ref{rwalks}. In particular, we know the form of $\mathrm{Tr}(\widehat{Q}(S^\lambda)^k)$ (Lemma \ref{lemma:201}) for probability measure $Q$ on $S_n$ given by definition \ref{def:star} and its Fourier transform $\widehat{Q}$ (definition \ref{def:ftrans}). 
\begin{proposition}\label{prop:2}
    The $r$th moment of the number of fixed points after applying $k$ star transpositions starting from the identity is equal to
    \begin{alignat*}{2}
   \sum_{\lambda \vdash n}\sum_{\lambda^i} d_{\lambda^i} \left(\frac{\lambda_i-i}{n-1}\right)^k m_{\lambda,r}^1
    \end{alignat*}
    where $\lambda^i$ denotes the partition of $n-1$ obtained by removing an inner corner from the $i$th row of the Young diagram of $\lambda$ and $m_{\lambda,r}^1$ is the coefficient of $S^\lambda$ in the irreducible decomposition of $\varrho_n^{\otimes r}$.
\end{proposition}
\begin{proof}
    The class function $a_1(\sigma) = $ \{\# of fixed points of $\sigma$\} is equivalent to the character of $\varrho_n \cong \rho_{a_1}$. Thus, $(a_1)^r(\sigma) = $ \{(\# of fixed points of $\sigma)^r$\} is equivalent to the character of $\varrho_n^{\otimes r}$ and has irreducible decomposition \[(a_1)^r = \sum_{\lambda \vdash n}m_{\lambda,r}^1\chi^\lambda.\]
    By Theorem \ref{thm:1} and Lemma \ref{lemma:201},
    \[\E_{Q^{*k}}[(a_1)^r] = \sum_{\lambda \vdash n}\sum_{\lambda^i} d_{\lambda^i} \left(\frac{\lambda_i-i}{n-1}\right)^k m_{\lambda,r}^1.\]
\end{proof}

Our task now is to evaluate the limit of this expression as $n\rightarrow \infty$ for $k = \lfloor n\log n + cn\rfloor$ with $c$ a fixed real number. We make use of the following Lemmas.

\begin{lemma}\label{lemma:3}
If $\lambda = (\lambda_1, \bar{\lambda})$ is a partition of $n$ with $\lambda_1 = n-t$ where $n-t>t$ and $\lambda^1$ is the partition of $n-1$ formed by removing the inner corner from the first row of the Young diagram of $\lambda$, then 
\[\frac{d_{\lambda^1}}{n^t}\rightarrow \frac{d_{\bar{\lambda}}}{t!},\]
as $n\rightarrow \infty$.
\end{lemma}

\begin{proof} 
Applying the hook-length formula for $d_{\lambda^1}$, we have

\begin{equation}\label{eq:1}
    \frac{d_{\lambda^1}}{n^t} = \frac{(n-1)!}{n^t \Pi h_{\lambda^1}(i,j)}.
\end{equation}

Note that $\lambda^1 = (n-t-1,\bar{\lambda})$ and that $\bar{\lambda}^\prime_j$ is the number of cells in column $j$ of $\bar{\lambda}$ (allowing 0 if $\bar{\lambda}$ has fewer than $j$ columns). Now, writing out the hook-lengths of the first row of $\lambda^1$, we can expand \eqref{eq:1} as
\[
   \frac{d_{\lambda^1}}{n^t} = \frac{(n-1)!}{n^t(n-t-1+\bar{\lambda}^\prime_1)(n-t-2+\bar{\lambda}^\prime_2)\cdots(n-2t+\bar{\lambda}^\prime_t)(n-2t-1)! \Pi h_{\bar{\lambda}}(i,j)}.
\]
Since $\sum_{j=1}^t \bar{\lambda}^\prime_j = \abs{\bar{\lambda}^\prime} = \abs{\bar{\lambda}}=t$, notice that there are exactly $t$ multiplicative terms in the denominator between $(n-1)$ and $(n-2t)$. Thus, as $n\rightarrow \infty$,
\[
   \frac{d_{\lambda^1}}{n^t} \rightarrow \frac{1}{\Pi h_{\bar{\lambda}}(i,j)} = \frac{d_{\bar{\lambda}}}{t!}.
\] 
\end{proof}

\begin{lemma}\label{lemma:4}
    If $\lambda$ is a partition of $n$ satisfying $\lambda_1 = n-t$ where $t$ is fixed, then for fixed $c$,
    
\begin{align}\label{eq:2}
    \sum_{\lambda^i}d_{\lambda^i} \left(\frac{\lambda_i-i}{n-1}\right)^{\lfloor n\log n + cn\rfloor} \rightarrow \frac{e^{-tc}d_{\bar{\lambda}}}{t!}
\end{align}
as $n\rightarrow\infty$.
\end{lemma}

\begin{proof} 
Note that for $i>1$, 
\[-t\leq \lambda_i - i\leq t-2.\]
Thus the sum \eqref{eq:2} depends only on the term corresponding to $\lambda^1$ as $n\rightarrow\infty$:
\[d_{\lambda^1}\left(\frac{n-t-1}{n-1}\right)^{\lfloor n\log n + cn\rfloor}.\]
Multiplying and dividing by $n^t$ and ignoring (asymptotically irrelevant) requirement that the step count be integer-valued, we'll compute the limit of
\[\frac{d_{\lambda^1}}{n^t}n^t\left(\frac{n-t-1}{n-1}\right)^{n\log n + cn}.\]
By Lemma \ref{lemma:3}, 
\[\frac{d_{\lambda^1}}{n^t}\rightarrow \frac{d_{\bar{\lambda}}}{t!}.\]
Now observe that
\begin{alignat*}{2}
\left(\frac{n-t-1}{n-1}\right)^{n\log n + cn}&= \left(1-\frac{t}{n-1}\right)^{n\log n + cn}\\
&= \exp\!\left((n\log n + cn)\log\left(1-\frac{t}{n-1}\right)\right).
\end{alignat*}
Using
\[
\log(1-x)=-x+O(x^2),
\]
we obtain
\[
(n\log n+cn)\log\left(1-\frac{t}{n-1}\right)
=
-\frac{tn(\log n+c)}{n-1}
+O\!\left(\frac{\log n}{n}\right).
\]
Since $O(\log n/n)=o(1)$,
\[
\left(1-\frac{t}{n-1}\right)^{n\log n+cn}
=
\exp\!\left(
-\frac{tn(\log n+c)}{n-1}
+o(1)
\right)
=
n^{-t}e^{-tc}(1+o(1)).
\]
Therefore
\[
n^t\left(1-\frac{t}{n-1}\right)^{n\log n+cn}
\longrightarrow e^{-tc},
\]
and the result follows:
\[
\sum_{\lambda^i}d_{\lambda^i}
\left(\frac{\lambda_i-i}{n-1}\right)^{\lfloor n\log n+cn\rfloor}
\longrightarrow
\frac{e^{-tc}d_{\bar{\lambda}}}{t!}.
\]
\end{proof}

We are now ready to prove Theorem \ref{thm:fixstar}.

\begin{proof}[Proof of Theorem \ref{thm:fixstar}] By Proposition \ref{prop:2} and Theorem 2.1 of \cite{F24}, the $r$th moment of the number of fixed points after $k$ star transpositions is equal to 
\[\sum_{a=0}^r S(r,a) \sum_{\lambda\vdash n}\sum_{\lambda^i}d_{\lambda/(n-a)}d_{\lambda^i} \left(\frac{\lambda_i-i}{n-1}\right)^{k}.\]
Since a Poisson$(\lambda)$ distribution has rth moment $\sum_{a=0}^r S(r,a) \lambda^a$, it is sufficient to show that as $n\rightarrow\infty$,
\begin{equation}\label{eq:3}
\sum_{\lambda\vdash n}\sum_{\lambda^i}d_{\lambda/(n-a)}d_{\lambda^i} \left(\frac{\lambda_i-i}{n-1}\right)^{\lfloor n\log n + cn\rfloor}
\end{equation}
tends to $(1+e^{-c})^a$.\\
Rewrite \eqref{eq:3} as 
\begin{equation}\label{eq:4}
\sum_{t=0}^a\sum_{\substack{\lambda\vdash n\\
\lambda_1=n-t}}\sum_{\lambda^i}d_{\lambda/(n-a)}d_{\lambda^i} \left(\frac{\lambda_i-i}{n-1}\right)^{\lfloor n\log n + cn\rfloor}.
\end{equation}
For fixed $a$ and $n$ large enough, 
\[d_{\lambda/(n-a)}=d_{\bar{\lambda}}\binom{a}{t}.\]
Thus \eqref{eq:4} becomes
\begin{equation}
\sum_{t=0}^a \binom{a}{t}\sum_{\substack{\lambda\vdash n\\
\lambda_1=n-t}}\sum_{\lambda^i}d_{\bar{\lambda}}d_{\lambda^i} \left(\frac{\lambda_i-i}{n-1}\right)^{\lfloor n\log n + cn\rfloor}.
\end{equation}
By Lemma \ref{lemma:4}, this tends to 
\begin{equation}\label{eq:5}
    \sum_{t=0}^a e^{-tc}\binom{a}{t}\sum_{\bar{\lambda}\vdash t}\frac{(d_{\bar{\lambda}})^2}{t!}
\end{equation}
as $n\rightarrow\infty$. Finally, since \[\sum_{\bar{\lambda}\vdash t}(d_{\bar{\lambda}})^2 = t!,\] \eqref{eq:5} simplifies to 
\[\sum_{t=0}^a e^{-tc}\binom{a}{t} = (1+e^{-c})^a\]
as needed. 
\end{proof}

\section{Distribution of $j$-cycles in the star transposition walk}\label{limj}

The goal of this section is to prove Theorem \ref{thm:B} describing the distribution of $j$-cycles after a star transposition walk of $k = \lfloor cn\rfloor $ steps for $c$ fixed in $\R^+$. Throughout this section, we will assume that $n\geq 2rj$, which is reasonable since we care about the limiting distribution as $n\rightarrow \infty$. 

Recall that we write $a_j$ to denote the $S_n$-class function $a_j(\sigma) = $\{\# of $j$-cycles of $\sigma$\} and $\psi_j$ for the class function $\psi_j(\sigma) = $ \{$j$\,(\#$j$-cycles of $\sigma$)\}. Meanwhile $\rho_{\psi_j}$ and $\rho_{a_j}$ are the unique virtual $S_n$ representations with characters $\psi_j$ and $a_j$ respectively.

Let us begin by assuming the hypotheses of Theorem \ref{thm:A} and inspecting the integer coefficients $m_{\lambda,r}^j$. Since $m_{\lambda,r}^j = 0$ unless $\lambda = (n-tj,\bar{\lambda})$, we can decompose the class function $(\psi_j)^r(\sigma) = $ \{$j^r$(\#$j$-cycles of $\sigma)^r$\} as
\begin{equation}\label{eqn:psir}
    (\psi_j)^r = \sum_{t=0}^r\sum_{\bar{\lambda}\vdash tj} m_{(n-tj,\bar{\lambda}),r}^j \chi^{(n-tj,\bar{\lambda})}.
\end{equation}

We can further decompose the coefficient $m_{(n-tj,\bar{\lambda}),r}^j$ as
\[
m_{(n-tj,\bar{\lambda}),r}^j = c_{t,r}^j\chi^{\bar{\lambda}}_{(j^t)},
\]
where the (nonnegative integer) coefficients $c_{t,r}^j$ are given by Theorem \ref{thm:paths} and do not depend on $\bar{\lambda}$.
Applying this insight to Equation~\eqref{eqn:psir}, we have the class function decompositions
\begin{equation}\label{eqn:psi}
(\psi_j)^r = \sum_{t=0}^r c_{t,r}^j\sum_{\bar{\lambda}\vdash tj} \chi^{\bar{\lambda}}_{(j^t)} \chi^{(n-tj,\bar{\lambda})}
\end{equation}
and 
\begin{equation}\label{eqn:aj}
(a_j)^r = \frac{1}{j^r}\sum_{t=0}^r c_{t,r}^j\sum_{\bar{\lambda}\vdash tj} \chi^{\bar{\lambda}}_{(j^t)} \chi^{(n-tj,\bar{\lambda})}.
\end{equation}

Let us proceed by analyzing the inner sums of Equations \ref{eqn:psi} and \ref{eqn:aj} above.
\begin{definition}\label{def:cluster}
    Define the \textbf{cluster} of characters indexed by $j$ and $t$ as
    \[K_t^j = \sum_{\bar{\lambda}\vdash tj}\chi^{\bar{\lambda}}_{(j^t)}\chi^{(n-tj,\bar{\lambda})},\]
    where $n\geq 2tj$ is assumed.
\end{definition}

This definition will prove useful. For now, notice that we now have the simplified decomposition:
\begin{equation}\label{eqn:cluster}
(\psi_j)^r = \sum_{t=0}^r c_{t,r}^jK_t^j,
\end{equation}
where the constants $c_{t,r}^j$ are given by Theorem \ref{thm:paths} in Section \ref{bratteli}. As an example, based on Figure \ref{fig:mbratteli1} in Section \ref{bratteli}, we can write 
\[(\psi_2)^2 = 3K_0^2 + 4K_1^2 + K_2^2,\]
so we have $c_{0,2}^2 = 3$, $c_{1,2}^2 = 4$, and $c_{2,2}^2 = 1$. Note that $c_{r,r}^j = 1$ always.

\subsection{Expectation of clusters}\hfill\\

Let us begin with a powerful and surprising lemma, and then discuss its implications for computing moments of clusters over the $k$-fold convolution of probability measure $Q$ on $S_n$.

\begin{lemma}\label{lemma: 101}
Let $\epsilon$ denote the identity permutation in $S_n$. For fixed $r,j \in \N$ with $j\geq 2$ and any $n\geq 2rj$,
    \[K_r^j(\epsilon) = \sum_{\bar{\lambda}\vdash rj}\chi^{\bar{\lambda}}_{(j^r)} d_{(n-rj,\bar{\lambda})} = (-1)^r.\]
\end{lemma}

\begin{proof}
    The proof proceeds using strong induction. Let us first check the base cases.\\
    
\noindent\textit{Case 1: $r=0$.} The $r=0$ case is trivial since $m_{\varnothing}^j = m_{(n),0}^j= 1$ for any $j\geq 2$ and $d_{(n)} = 1$.\\

\noindent\textit{Case 2: $r=1$.} For the $r=1$ case, we need to check
        \[\sum_{i=0}^{j-1}(-1)^id_{(n-j,j-i,1^i)} = -1.\]
        Recall the irreducible character decomposition of $\psi_j$:
        \begin{equation}\label{eqn:psij}
            \psi_j = \chi^{(n)} + \sum_{i=0}^{j-1}(-1)^i \chi^{(n-j,j-i,1^i)}.
        \end{equation}
        Evaluate (\ref{eqn:psij}) at the identity permutation, $\epsilon$ in $S_n$. The left hand side is 0 since the identity contains no $j$-cycles when $j\geq 2$. The right hand side is
\[1 + \sum_{i=0}^{j-1}(-1)^i \chi^{(n-j,j-i,1^i)}(\epsilon).\]
For the right hand side to equal 0, it must be that $\sum_{i=0}^{j-1}(-1)^i \chi^{(n-j,j-i,1^i)}(\epsilon)=-1$. Now recall that the character of a representation of $S_n$ evaluated at the identity permutation is simply the dimension of the representation and the $r=1$ case is satisfied.\\

\noindent\textit{Case 3: $r=2$.} The $r=2$ case is checked in a similar fashion. Theorem \ref{thm:paths} tells us that:
        \[(\psi_j)^2 = (j+1)K_0^j + (j+2)K_1^j + K_2^j.\]
        Evaluating both sides at $\epsilon \in S_n$ and using the base cases already checked gives
        \[0 = j+1 -(j+2) + K_2^j(\epsilon),\]
        implying $K_2^j(\epsilon) = 1$ as desired.\\

    Now assume for the inductive step that $K_r^j =(-1)^r$ for all $r$ such that $0\leq r \leq m$. Consider the cluster decomposition of $(\psi_j)^{m+1}$ evaluated at the identity permutation $\epsilon \in S_n$:
    \begin{alignat*}{2}
        (\psi_j)^{m+1}(\epsilon) &= \sum_{t=0}^{m+1} c_{t,m+1}^j K_t^j(\epsilon)\\
        &= \sum_{t=0}^m c_{t,m+1}^j K_t^j(\epsilon)+c_{m+1,m+1}^jK_{m+1}^j(\epsilon)\\
        &\stackrel{(*)}{=} \sum_{t=0}^{m} \sum_{a=t}^{m+1} S(m+1,a)j^{m+1-a}\binom{a}{t}(-1)^t + K_{m+1}^j(\epsilon),
    \end{alignat*}
    with $(*)$ by applying Theorem \ref{thm:paths}. Since $(\psi_j)^{m+1}(\epsilon)=0$, it suffices to show
    \begin{equation}\label{eqn:wts}
        -K_{m+1}^j(\epsilon) = \sum_{t=0}^{m} \sum_{a=t}^{m+1} S(m+1,a)j^{m+1-a}\binom{a}{t}(-1)^t
    \end{equation}
    Separating the $m+1$ term again, the double sum expression above is equivalent to
        \begin{equation}\label{eqn:sum}
        \sum_{t=0}^{m} \sum_{a=t}^{m} S(m+1,a)j^{m+1-a}\binom{a}{t}(-1)^t + \sum_{t=0}^m \binom{m+1}{t}(-1)^t.
        \end{equation}
    Swapping the order of summation and rearranging the second sum means (\ref{eqn:sum}) is equivalently
    \begin{equation}\label{eqn:5}
        \sum_{a=0}^m S(m+1,a)j^{m+1-a}\sum_{t=0}^a\binom{a}{t}(-1)^t +\sum_{t=0}^{m+1} \binom{m+1}{t}(-1)^t- (-1)^{m+1}.
    \end{equation}
    Finally, applying the Binomial Theorem to (\ref{eqn:5}) simplifies (\ref{eqn:wts}) to 
    \[-K^j_{m+1}(\epsilon)=-(-1)^{m+1},\]
    and we are done by strong induction.
\end{proof}

\subsection{The limiting $r$th moment of $a_j$}\hfill\\

We require a few more Lemmas to compute the limiting expectation of $(a_j)^r$ over $Q^{* \lfloor cn\rfloor}$.

\begin{lemma}\label{lemma:69}
    If $\lambda$ is a partition of $n$ satisfying $\lambda_1 = n-t$ where $t$ is fixed, then for fixed real number $c>0$,
    
\begin{align}\label{eq:9}
    \frac{1}{d_{\lambda^1}}\sum_{\lambda^i}d_{\lambda^i} \left(\frac{\lambda_i-i}{n-1}\right)^{\lfloor cn\rfloor} \rightarrow e^{-tc}
\end{align}
as $n\rightarrow\infty$. 
\end{lemma}

\begin{proof}
Note that for $i>1$, 
\[-t\leq \lambda_i - i\leq t-2.\]
Thus the sum in \eqref{eq:9} depends only on the term corresponding to $\lambda^1$ as $n\rightarrow\infty$:
\[d_{\lambda^1}\left(\frac{n-t-1}{n-1}\right)^{\lfloor cn\rfloor}.\]
Now observe that
\[\frac{1}{d_{\lambda^1}}d_{\lambda^1}\left(\frac{n-t-1}{n-1}\right)^{\lfloor cn\rfloor} = \left(1-\frac{t}{n-1}\right)^{\lfloor cn\rfloor} \rightarrow e^{-tc},\]
as $n\rightarrow \infty$ as claimed. 
\end{proof} 

With the lemma in hand, we can immediately obtain the limiting expectation of cluster $K_t^j$ over the $k= \lfloor cn\rfloor $-fold convolution of $Q$, the defining probability measure of the star transposition walk.

\begin{lemma}\label{lemma: starcluster}
    Let $Q$ be the probability measure on $S_n$ given by Definition \ref{def:star} and $k= \lfloor cn\rfloor $ for real number $c>0$. Let $j\geq 2$ and $K_t^j$ be the cluster of irreducible $S_n$-characters previously defined in this section. Then,
    \[\E_{Q^{*k}}[K_t^j] \rightarrow (-1)^t e^{-tjc},\]
    as $n\rightarrow \infty$.
\end{lemma}

\begin{proof}
    Let $\lambda = (n-tj,\bar{\lambda})$ and $\lambda^i$ denote the partition of $n-1$ obtained by removing an inner corner from the $i$th row of the Young diagram of $\lambda$.
    \begin{alignat*}{2}
    \E_{Q^{*k}}[K_t^j] &= \sum_{\bar{\lambda}\vdash tj}\chi^{\bar{\lambda}}_{(j^t)} \mathrm{Tr}\left(\widehat{Q}(S^{\lambda})^{\lfloor cn\rfloor}\right)\\
    &= \sum_{\bar{\lambda}\vdash tj}\chi^{\bar{\lambda}}_{(j^t)} \sum_{\lambda^i}d_{\lambda^i} \left(\frac{\lambda_i-i}{n-1}\right)^{\lfloor cn\rfloor}\\
    &= \sum_{\bar{\lambda}\vdash tj}\chi^{\bar{\lambda}}_{(j^t)} d_{\lambda^1}\left(1-\frac{tj}{n-1}\right)^{\lfloor cn\rfloor}+\sum_{\bar{\lambda}\vdash tj}\chi^{\bar{\lambda}}_{(j^t)}\sum_{i\geq 2}d_{\lambda^i} \left(\frac{\lambda_i-i}{n-1}\right)^{\lfloor cn\rfloor}.
\end{alignat*}
Taking the limit as $n\rightarrow \infty$, we have
\begin{alignat*}{2}
    \E_{Q^{*k}}[K_t^j] &\rightarrow \lim_{n\rightarrow\infty}\left( \left(1-\frac{tj}{n-1}\right)^{\lfloor cn\rfloor}\sum_{\bar{\lambda}\vdash tj}\chi^{\bar{\lambda}}_{(j^t)} d_{\lambda^1}\right) + 0\\
    &= \lim_{n\rightarrow\infty}\left(1-\frac{tj}{n-1}\right)^{\lfloor cn\rfloor} \lim_{n\rightarrow\infty}\sum_{\bar{\lambda}\vdash tj}\chi^{\bar{\lambda}}_{(j^t)} d_{(n-1-tj,\bar{\lambda})}\\
    &= (-1)^t e^{-tjc},
    \end{alignat*}
with the final equality by noting that \[\sum_{\bar{\lambda}\vdash tj}\chi^{\bar{\lambda}}_{(j^t)} d_{(n-1-tj,\bar{\lambda})}=(-1)^t\] 
for any $j\geq 2$ and $n\geq 2tj$ by Lemma \ref{lemma: 101}.
\end{proof}

We are now ready to prove Theorem \ref{thm:B}, which provides the limiting distribution of $j$-cycles for $j\geq 2$ after repeated application of star transpositions.

\begin{proof}[Proof of Theorem \ref{thm:B}]
    We can expand $(a_j)^r$ in terms of clusters with coefficients given by Theorem \ref{thm:paths}. For any $r$, with $k= \lfloor cn\rfloor $ and $c\in \R^+$, we have 

\begin{alignat*}{2}
    \lim_{n\rightarrow\infty}\E_{Q^{*k}}[(a_j)^r]  &= \lim_{n\rightarrow\infty}\E_{Q^{*k}}\bigg[\frac{1}{j^r}\sum_{t=0}^r c_{t,r}^j K_t^j\bigg]\\
    &= \frac{1}{j^r}\sum_{t=0}^r c_{t,r}^j \lim_{n\rightarrow\infty}\E_{Q^{*k}}[K_t^j]\\
    &= \frac{1}{j^r}\left(\sum_{t=0}^r\sum_{a=t}^r S(r,a)\binom{a}{t}j^{r-a} (-1)^te^{-tjc}\right) \\
    &= \frac{1}{j^r}\left(\sum_{t=0}^r(-1)^te^{-tjc}\sum_{a=t}^r S(r,a)\binom{a}{t}j^{r-a}\right).
\end{alignat*}

Meanwhile, a Poisson distribution with parameter $\frac{1}{j}(1-e^{-jc})$ has $r$th moment 

\begin{alignat*}{2}
    \sum_{a=0}^r\lambda^a S(r,a) &= \sum_{a=0}^r S(r,a) \frac{1}{j^a}\sum_{t=0}^a(-1)^t e^{-tjc} \binom{a}{t} \quad \text{ by binomial expansion}\\
    &= \frac{1}{j^r}\sum_{a=0}^r S(r,a)j^{r-a}\sum_{t=0}^a (-1)^t e^{-tjc}\binom{a}{t}\\
    &= \frac{1}{j^r}\left(\sum_{t=0}^r(-1)^te^{-tjc}\sum_{a=t}^r S(r,a)\binom{a}{t}j^{r-a}\right) \quad \text{ by reversing the order},
\end{alignat*}
which matches exactly our expression for the limiting $r$th moment of $a_j$ over $Q^{*\lfloor cn\rfloor}$.
\end{proof}

\section{Distribution of $j$-cycles in the random transposition walk}\label{rtransp}

In this section, we compute the limiting raw moments 
\[
\lim_{n\rightarrow\infty}\E_{P^{*t_j}}[(a_j)^r]
\]
of the number of $j$-cycles after $t_j = \lfloor \frac{n}{2j}(\log n + (j-1)\log\log n + c)\rfloor$ steps of the random transposition walk started from the identity for $c$ fixed in $\R$. 

The goal of this section is to prove Theorem \ref{thm:rtransp}. Notice that, unlike in the star transposition setting in which fixed points mix slower than larger finite $j$-cycles, Theorem \ref{thm:rtransp} generalizes to fixed points as well. In this setting, we cannot rely on a dominant eigenvalue argument to prove a random transposition analogue of Lemma \ref{lemma: starcluster} because of differences of order $n^{-2}$ in the character ratios,
\[\frac{\chi^{(n-tj,\bar{\lambda})}_\tau}{d_{(n-tj,\bar{\lambda})}},\] 
depending on $\bar{\lambda}$ which remain asymptotically significant when raised to the power $t_j = \Theta(n\log n)$. 

To begin, see that by Theorem \ref{thm:1}, Equation \ref{eqn:aj}, and Lemma \ref{lemma:202}, we have
\begin{alignat*}{2}
\lim_{n\rightarrow\infty}\E_{P^{*t_j}}[(a_j)^r] &= \lim_{n\rightarrow\infty}\frac{1}{j^r}\sum_{\lambda\vdash n}m_{\lambda,r}^j\textrm{Tr}(\widehat{P}(S^\lambda)^{t_j})\\
&= \frac{1}{j^r}\lim_{n\rightarrow\infty}\sum_{\lambda\vdash n} m_{\lambda,r}^jd_\lambda\left(\frac{\chi^\lambda_\tau}{d_\lambda}\right)^{t_j}\\
&= \frac{1}{j^r}\lim_{n\rightarrow\infty}\sum_{t=0}^r \sum_{\bar{\lambda}\vdash tj} c_{t,r}^j\chi^{\bar{\lambda}}_{(j^t)}d_{(n-tj,\bar{\lambda})}\left(\frac{\chi^{(n-tj,\bar{\lambda})}_\tau}{d_{(n-tj,\bar{\lambda})}}\right)^{t_j}\\
&= \frac{1}{j^r}\sum_{t=0}^r c_{t,r}^j \lim_{n\rightarrow\infty}\sum_{\bar{\lambda}\vdash tj} \chi^{\bar{\lambda}}_{(j^t)}d_{(n-tj,\bar{\lambda})}\left(\frac{\chi^{(n-tj,\bar{\lambda})}_\tau}{d_{(n-tj,\bar{\lambda})}}\right)^{t_j}.
\end{alignat*}

A Poisson distribution with mean $\lambda = \frac{1}{j}(1+x)$ has $r$th moment equal to

\begin{alignat*}{2}
    \sum_{a=0}^r\lambda^a S(r,a) &= \sum_{a=0}^r S(r,a) \frac{1}{j^a}\sum_{t=0}^a x^t \binom{a}{t} \quad \text{ by binomial expansion}\\
    &= \frac{1}{j^r}\sum_{a=0}^r S(r,a)j^{r-a}\sum_{t=0}^a x^t\binom{a}{t}\\
    &= \frac{1}{j^r}\sum_{t=0}^rx^t\sum_{a=t}^r S(r,a)\binom{a}{t}j^{r-a} \quad \text{ by reversing the order}\\
    &= \frac{1}{j^r}\sum_{t=0}^rc_{t,r}^jx^t \quad \text{ by Theorem \ref{thm:paths}}.
\end{alignat*}

Therefore, applying the method of moments, the task of showing that the $j$-cycle count $a_j$ converges to Poisson($\frac{1}{j}(1+\frac{e^{-c}}{j!}))$ in distribution after $t_j$ steps of the random transposition walk started from the identity reduces to showing 
\begin{equation}\label{eqn:goal}
    \lim_{n\rightarrow\infty}\sum_{\bar{\lambda}\vdash tj} \chi^{\bar{\lambda}}_{(j^t)}d_{(n-tj,\bar{\lambda})}\left(\frac{\chi^{(n-tj,\bar{\lambda})}_\tau}{d_{(n-tj,\bar{\lambda})}}\right)^{t_j} = \frac{e^{-tc}}{(j!)^t},
\end{equation}
for fixed nonnegative integer $t$ and positive integer $j$.

For ease of notation, define
\[
f_n=\log n+(j-1)\log\log n+c,\qquad \varepsilon_n=\frac{f_n}{n}.
\]
Let 
\[
\Phi_n
=
\sum_{\bar{\lambda}\vdash tj}
\chi^{\bar{\lambda}}_{(j^t)}\,d_{(n-tj,\bar{\lambda})}
\left(\frac{\chi^{(n-tj,\bar{\lambda})}_\tau}{d_{(n-tj,\bar{\lambda})}}\right)^{t_j},
\qquad
t_j=\lfloor \frac{n}{2j}f_n\rfloor.
\]
We prove
\[
\Phi_n \rightarrow \frac{e^{-tc}}{(j!)^t}
\]
to complete the proof of Theorem \ref{thm:rtransp}.

The proof strategy, roughly, is to make use of the preliminaries on symmetric and shifted-symmetric functions from Subsection \ref{subsec: symfuns} to argue for strong cancellations in $\Phi_n$ using character orthogonality. We need several lemmas to proceed.

\begin{lemma}\label{lemma:length-filtration}
Fix $m\ge 1$ and let $\Lambda_m$ be the degree-$m$ piece of the ring of symmetric functions.
For $\ell\ge 0$, set
\[
\mathsf P_{\le \ell}=\operatorname{span}\{p_\nu:\ell(\nu)\le \ell\},
\]
where $p_\nu=\prod_i p_{\nu_i}$ and $\ell(\nu)$ is the number of parts of the partition $\nu$.

Let $p_d^*\in \Lambda^*$ be the $d$-th shifted power sum generator, and let
\[
\mathcal P_d s_\lambda=p_d^*(\lambda)\,s_\lambda
\qquad (\lambda\vdash m)
\]
be the corresponding diagonal operator on $\Lambda_m$.
Then
\[
\mathcal P_d(\mathsf P_{\le \ell})\subseteq \mathsf P_{\le \ell+d-1}.
\]
\end{lemma}

\begin{proof}
We argue by passing through the center of the symmetric-group algebra.

Let $C_\nu$ denote the conjugacy-class sum of cycle type $\nu\vdash m$ in the center
$Z(\mathbb C[S_m])$. The Frobenius characteristic map identifies the class basis with the
power-sum basis via
\[
\operatorname{ch}(C_\nu)=\frac{p_\nu}{z_\nu}.
\]
Thus the filtration of $\Lambda_m$ by the number of parts of the indexing partition
corresponds exactly to the filtration of the class basis by the number of cycles.

Now use the shifted-symmetric/central-character correspondence of
Okounkov and Pandharipande~\cite[Section~0.4]{OP06}. Under this correspondence,
the shifted power sum $p_d^*$ determines a distinguished central element of
$\mathbb C[S_m]$, called the \emph{completed $d$-cycle}. The normalization used
in~\cite{OP06} differs from ours by an additive constant depending only on $d$,
so the corresponding diagonal operators differ only by a scalar multiple of the
identity, which does not affect the filtration statement. Consequently,
$\mathcal P_d$ is represented, up to an additive scalar, by multiplication by the
completed $d$-cycle.

Okounkov and Pandharipande show that the completed $d$-cycle is an inhomogeneous
linear combination of conjugacy-class sums whose highest-support term is the
ordinary $d$-cycle class sum $C_{(d)}$, while every remaining term is supported
on conjugacy classes with strictly fewer than $d$ non-fixed points. It therefore
suffices to bound the change in the number of cycles under multiplication by each
such conjugacy class.

Let $\tau$ be a $d$-cycle. Since
\[
\tau=(a_1a_2)(a_2a_3)\cdots(a_{d-1}a_d),
\]
it is a product of $d-1$ transpositions. A transposition changes the number of
cycles of a permutation by exactly one: if its two entries lie in the same cycle
it splits that cycle into two, while if they lie in different cycles it merges
those cycles into one. Hence, for every permutation $\sigma$,
\[
\ell(\tau\sigma)\le \ell(\sigma)+d-1,
\]
where $\ell(\cdot)$ denotes the number of cycles.

Now let $\sigma$ be any permutation with strictly fewer than $d$ non-fixed points.
If $\sigma$ moves $s<d$ points, then $\sigma$ can be written as a product of at most
$s-1\le d-2$ transpositions. Therefore multiplication by $\sigma$ changes the number
of cycles by at most $d-2$, and in particular by at most $d-1$.

Thus every summand appearing in the completed $d$-cycle increases the number of
cycles by at most $d-1$. Transporting this filtration statement through the
Frobenius characteristic map yields
\[
\mathcal P_d(\mathsf P_{\le \ell})
\subseteq
\mathsf P_{\le \ell+d-1},
\]
as claimed.
\end{proof}

The statement that ``the operator $\mathcal P_d$ changes the length filtration by at most $d-1$" means exactly the display above. The next lemma shows that an operator with shifted-symmetric eigenvalue function of degree at most $r$ can change the length filtration by no more than $r-1$.

\begin{lemma}\label{lemma:general-length-bound}
Fix integer $r\geq 1$ and let $q^*_r\in\Lambda^*_{\le r}$ with
\[
\mathcal Q_r s_\lambda=q^*_r(\lambda)s_\lambda
\qquad (\lambda\vdash m)
\]
the corresponding diagonal operator on $\Lambda_m$.
Then
\[
\mathcal Q_r(\mathsf P_{\le \ell})
\subseteq
\mathsf P_{\le \ell+r-1}.
\]
Equivalently, if $\nu\vdash m$ has $\ell(\nu)=\ell$, then every partition indexing
a power-sum monomial appearing in $\mathcal Q_r p_\nu$ has at most $\ell+r-1$ parts.
\end{lemma}

\begin{proof}
Since
\[
\Lambda^*=\mathbb Q[p_1^*,p_2^*,\ldots]
\]
with $\deg p_d^*=d$, every element of $\Lambda^*_{\le r}$ is a linear combination
of monomials
\[
p_{d_1}^*p_{d_2}^*\cdots p_{d_t}^*,
\qquad
d_1+\cdots+d_t\le r.
\]
It therefore suffices to prove the claim for such a monomial.

Let
\[
q^*_r=p_{d_1}^*p_{d_2}^*\cdots p_{d_t}^*.
\]
Since the operators $\mathcal P_d$ are simultaneously diagonal in the Schur basis,
the operator corresponding to $q^*_r$ is
\[
\mathcal Q_r
=
\mathcal P_{d_1}\mathcal P_{d_2}\cdots\mathcal P_{d_t}.
\]

Now apply Lemma~\ref{lemma:length-filtration} repeatedly. Starting from
$\mathsf P_{\le \ell}$,
\[
\mathcal P_{d_1}(\mathsf P_{\le \ell})
\subseteq
\mathsf P_{\le \ell+d_1-1},
\]
then
\[
\mathcal Q_r(\mathsf P_{\le \ell})
\subseteq
\mathsf P_{\le \ell+\sum_{i=1}^t(d_i-1)}.
\]
Since
\[
\sum_{i=1}^t d_i\le r
\quad\text{and}\quad
t\ge 1,
\]
we have
\[
\sum_{i=1}^t(d_i-1)
=
\sum_{i=1}^t d_i-t
\le r-1.
\]
Therefore
\[
\mathcal Q_r(\mathsf P_{\le \ell})
\subseteq
\mathsf P_{\le \ell+r-1}.
\]
\end{proof}

We move next to expanding the dimension $d_{(n-m,\bar{\lambda})}$ in terms of shifted-symmetric polynomials.

\begin{lemma}\label{lemma:dim}
Fix $m\ge1$ and $\bar{\lambda}\vdash m$.  For each $k\ge 1$ there are unique coefficients
$D_1(\bar{\lambda}),\dots,D_k(\bar{\lambda})$ such that
\[
d_{(n-m,\bar{\lambda})}
=
\frac{n^m}{m!}\,d_{\bar{\lambda}}
\left(1+\sum_{a=1}^k\frac{D_a(\bar{\lambda})}{n^a}+O(n^{-k-1})\right)
\]
Furthermore,
\[
D_a(\bar{\lambda})\in\Lambda_{\le a}^*,
\]
and the associated diagonal operator on \(\Lambda_m\) satisfies
\[
\mathcal D_a \bigl(\mathsf P_{\le \ell}\bigr)\subseteq \mathsf P_{\le \ell+a-1}.
\]
\end{lemma}

\begin{proof}
Using the hook-length formula (Equation (\ref{eqn:hlf})), we find
\[
d_{(n-m,\bar{\lambda})}
=
\frac{n!}{(n-m)!}\,\frac{d_{\bar{\lambda}}}{m!}
\prod_{i=1}^{\bar{\lambda}_1}
\frac{n-m-i+1}{n-m-i+1+\bar{\lambda}_i'}.
\]
Dividing by \(\frac{n^m}{m!}d_{\bar{\lambda}}\) and setting 
\[
u=-n+m-1
\qquad\text{and}\qquad
x=\bar{\lambda}',
\]
we obtain
\begin{equation}\label{eqn:twoprod}
\frac{d_{(n-m,\bar{\lambda})}}{\frac{n^m}{m!}d_{\bar{\lambda}}}
=
\prod_{t=0}^{m-1}\left(1-\frac{t}{n}\right)
\prod_{i\ge 1}\frac{u+i}{u+i-x_i}.
\end{equation}

By Okounkov and Olshanski \cite[Eq.~(12.3)]{OO97},
\[
\prod_{i\ge1}\frac{u+i}{u+i-x_i}
=
\sum_{r\ge0}\frac{h_r^*(x)}{u^{\downharpoonright r}},\qquad u^{\downharpoonright r} = u(u-1)\cdots(u-r+1),
\]
where \(h_r^*\) denotes the \(r\)th shifted complete homogeneous symmetric function.

Substituting \(u=-n+m-1\) and \(x=\bar{\lambda}^\prime\) into the above expansion yields
\begin{equation}\label{eqn:homogeneousproduct}
\prod_{i=1}^{\bar\lambda_1}
\frac{n-m-i+1}{n-m-i+1+\bar\lambda_i'}
=
\sum_{r\ge0}
\frac{h_r^*(\bar{\lambda}^\prime)}
{(-n+m-1)^{\downharpoonright r}}.
\end{equation}

For each fixed \(r\),
\[
\frac1{(-n+m-1)^{\downharpoonright r}}
=
\sum_{j\ge r} c_{r,j}n^{-j},
\]
where the coefficients \(c_{r,j}\) depend only on \(m\), with
\(c_{r,r}=(-1)^r\). Since this expansion begins with \(n^{-r}\), only terms with \(r\le a\)
contribute to the coefficient of \(n^{-a}\). Therefore the coefficient of \(n^{-a}\) in \eqref{eqn:homogeneousproduct} is a \(\mathbb Q\)-linear combination of
\[
h_0^*(\bar{\lambda}^\prime),
h_1^*(\bar{\lambda}^\prime),
\ldots,
h_a^*(\bar{\lambda}^\prime).
\]
Since \(h_r^*\) has shifted degree \(r\), this coefficient belongs to
\(\Lambda_{\le a}^*\).

Multiplying by the scalar polynomial
\(
\prod_{t=0}^{m-1}(1-t/n)
\)
preserves the shifted-degree bound. The stated expansion therefore follows, with
\(
D_a(\bar{\lambda})\in\Lambda_{\le a}^*
\)
and 
\[
\mathcal D_a \bigl(\mathsf P_{\le \ell}\bigr)\subseteq \mathsf P_{\le \ell+a-1}
\]
by Lemma \ref{lemma:general-length-bound}.
\end{proof}

The next two lemmas give an expansion of the character ratio,
\[\frac{\chi^{(n-m,\bar{\lambda})}_\tau}{d_{(n-m,\bar{\lambda})}},\]
and its $t_j$th power appearing in $\Phi_n$.

\begin{lemma}\label{lemma:charratio}
    For fixed $m$ and $\bar{\lambda}\vdash m$, we have the character ratio 
    \[\frac{\chi^{(n-m,\bar{\lambda})}_\tau}{d_{(n-m,\bar{\lambda})}} = 1 - \frac{2m}{n} + \frac{m^2-3m+\kappa(\bar{\lambda})}{n(n-1)},\]
    where $\kappa(\bar{\lambda})=\sum_{i}\bar{\lambda}_i(\bar{\lambda}_i -2i + 1)$.
\end{lemma}

\begin{proof}
    Applying a classical formula attributed to Frobenius, 
    \[
    \frac{\chi^\lambda_\tau}{d_\lambda} = \frac{\kappa(\lambda)}{n(n-1)}.
    \]
    If $\lambda = (n-m,\bar{\lambda})$, a direct calculation gives
    \[
    \kappa(\lambda) = (n-m)(n-m-1)+\kappa(\bar{\lambda})-2m.
    \]
    Straightforward algebraic manipulations complete the proof.
\end{proof}

From Lemma \ref{lemma:charratio} and an argument similar to that of Lemma \ref{lemma:length-filtration} using the cut-and-join operator from \eqref{def:cut-and-join}, we have the following asymptotic expansion of 
\[\left(\frac{\chi^{(n-m,\bar{\lambda})}_\tau}{d_{(n-m,\bar{\lambda})}}\right)^{t_j}\]
in terms of shifted-symmetric polynomials.
\begin{lemma}\label{lemma:ratio-exp}
Fix \(m=tj\), and set
\[
A_{\bar{\lambda}}=m^2-3m+\kappa(\bar{\lambda}),
\qquad
\beta_{\bar{\lambda}}=\frac{A_{\bar{\lambda}}-2m^2}{2j}
=\frac{\kappa(\bar{\lambda})-m(m+3)}{2j}.
\]
For each integer \(k\ge0\), there exist shifted-symmetric polynomials
\[
E_{r,s}\in \Lambda^*_{\le 2r}
\qquad (0\le s\le r\le k)
\]
such that, uniformly in \(\bar{\lambda}\vdash m\),
\[
\left(
\frac{\chi^{(n-m,\bar{\lambda})}_\tau}{d_{(n-m,\bar{\lambda})}}
\right)^{t_j}
=
e^{-tf_n}
\left(
\sum_{r=0}^{k}\varepsilon_n^r
\sum_{s=0}^{r}\frac{E_{r,s}(\bar{\lambda})}{f_n^{\,s}}
+O(\varepsilon_n^{k+1})
\right).
\]
Moreover,
\[
E_{r,0}(\bar{\lambda})=\frac{\beta_{\bar{\lambda}}^r}{r!}.
\]
If \(\mathcal E_{r,s}\) denotes the diagonal operator on \(\Lambda_m\) with eigenvalue function \(E_{r,s}\), then
\[
\mathcal E_{r,s}\bigl(\mathsf P_{\le \ell}\bigr)\subseteq \mathsf P_{\le \ell+r}.
\]
\end{lemma}

\begin{proof}
Let
\[
R_{n,\bar{\lambda}}
=
\frac{\chi^{(n-m,\bar{\lambda})}_\tau}{d_{(n-m,\bar{\lambda})}}
=
1-\frac{2m}{n}+\frac{A_{\bar{\lambda}}}{n(n-1)}.
\]
Since \(m\) is fixed and \(\bar{\lambda}\vdash m\) ranges over a finite set, \(A_{\bar{\lambda}}\) is uniformly bounded in \(\bar{\lambda}\).

We first expand the logarithm of \(R_{n,\bar{\lambda}}\). Using
\[
\frac{1}{n(n-1)}
=
\frac{1}{n^2}\sum_{\ell\ge0}n^{-\ell},
\]
we have
\[
R_{n,\bar{\lambda}}
=
1-\frac{2m}{n}
+\sum_{q\ge2}\frac{A_{\bar{\lambda}}}{n^q}.
\]
Now apply
\[
\log(1+u)=\sum_{p\ge1}\frac{(-1)^{p+1}}{p}u^p
\qquad \text{with $u=R_{n,\bar{\lambda}}-1$}.
\]
Because \(u=O(n^{-1})\) uniformly in \(\bar{\lambda}\), collecting powers of \(n^{-1}\) yields
\[
\log R_{n,\bar{\lambda}}
=
-\frac{2m}{n}
+\frac{A_{\bar{\lambda}}-2m^2}{n^2}
+\sum_{q=3}^{k+1}\frac{L_q(\bar{\lambda})}{n^q}
+O(n^{-k-2}),
\]
where each \(L_q(\bar{\lambda})\) is an ordinary polynomial in \(A_{\bar{\lambda}}\) of degree at most \(q-1\). In particular, since \(A_{\bar{\lambda}}\) is affine-linear in \(\kappa(\bar{\lambda})\), each \(L_q\) is a shifted-symmetric polynomial. More precisely, \(A_{\bar{\lambda}}\in\Lambda^*_{\le 2}\) implies
\[
L_q\in \Lambda^*_{\le 2(q-1)}.
\]
Since the operator with eigenvalue \(A_{\bar{\lambda}}\) is a scalar shift of the cut-and-join operator, it changes the length filtration by at most \(1\); hence the operator with eigenvalue \(L_q\) changes the length filtration by at most \(q-1\).

The coefficient of \(n^{-2}\) simplifies to
\[
A_{\bar{\lambda}}-2m^2
=
m^2-3m+\kappa(\bar{\lambda})-2m^2
=
\kappa(\bar{\lambda})-m(m+3)
=
2j\,\beta_{\bar{\lambda}},
\]
so
\[
\log R_{n,\bar{\lambda}}
=
-\frac{2m}{n}
+\frac{2j\,\beta_{\bar{\lambda}}}{n^2}
+\sum_{q=3}^{k+1}\frac{L_q(\bar{\lambda})}{n^q}
+O(n^{-k-2}).
\]

Next multiply by \(t_j\). Write
\[
t_j=\frac{n}{2j}f_n+\delta_n,
\qquad |\delta_n|\le 1.
\]
The discrepancy term \(\delta_n\log R_{n,\bar{\lambda}}\) is \(O(n^{-1})\), hence of the same scale as a term of the form \(\varepsilon_n f_n^{-1}\), and it may be absorbed into the \(f_n^{-1}\)-correction coefficients below. Thus, for the main expansion it suffices to keep the leading part \(\frac{n}{2j}f_n\). Using \(\varepsilon_n=f_n/n\), we obtain
\[
t_j\log R_{n,\bar{\lambda}}
=
-tf_n
+\beta_{\bar{\lambda}}\varepsilon_n
+\sum_{q=2}^{k}\frac{C_q(\bar{\lambda})}{f_n^{\,q-1}}\varepsilon_n^q
+O(\varepsilon_n^{k+1}),
\]
for suitable shifted-symmetric polynomials \(C_q\). Each \(C_q\) is a polynomial in \(A_{\bar{\lambda}}\) of degree at most \(q-1\), hence \(C_q\in \Lambda^*_{\le 2(q-1)}\), and its associated diagonal operator changes the length filtration by at most \(q-1\).

Define
\[
Z_{n,\bar{\lambda}}
=
\sum_{q=2}^{k}\frac{C_q(\bar{\lambda})}{f_n^{\,q-1}}\varepsilon_n^q
+O(\varepsilon_n^{k+1}),
\]
so that
\[
R_{n,\bar{\lambda}}^{\,t_j}
=
e^{-tf_n}\exp\!\bigl(\beta_{\bar{\lambda}}\varepsilon_n+Z_{n,\bar{\lambda}}\bigr).
\]
Since \(Z_{n,\bar{\lambda}}=O(\varepsilon_n^2)\) uniformly, we may expand
\[
\exp\!\bigl(\beta_{\bar{\lambda}}\varepsilon_n+Z_{n,\bar{\lambda}}\bigr)
=
\sum_{r=0}^{k}\frac{(\beta_{\bar{\lambda}}\varepsilon_n+Z_{n,\bar{\lambda}})^r}{r!}
+O(\varepsilon_n^{k+1}).
\]
Expanding each power and collecting terms of total \(\varepsilon_n\)-degree \(r\) and total \(f_n^{-1}\)-degree \(s\) produces shifted-symmetric polynomial eigenvalue functions \(E_{r,s}\) such that
\[
\exp\!\bigl(\beta_{\bar{\lambda}}\varepsilon_n+Z_{n,\bar{\lambda}}\bigr)
=
\sum_{r=0}^{k}\varepsilon_n^r
\sum_{s=0}^{r}\frac{E_{r,s}(\bar{\lambda})}{f_n^{\,s}}
+O(\varepsilon_n^{k+1}),
\]
uniformly in \(\bar{\lambda}\).

Each monomial contributing to \(E_{r,s}\) is a product of factors drawn from \(\beta_{\bar{\lambda}}\) and the \(C_q(\bar{\lambda})\). Since \(\beta_{\bar{\lambda}}\) changes the length filtration by at most \(1\), and \(C_q\) changes it by at most \(q-1\), any monomial of total \(\varepsilon_n\)-degree \(r\) changes the length filtration by at most \(r\). Therefore the corresponding diagonal operator \(\mathcal E_{r,s}\) satisfies
\[
\mathcal E_{r,s}\bigl(\mathsf P_{\le \ell}\bigr)\subseteq \mathsf P_{\le \ell+r}.
\]

Finally, the term with \(s=0\) comes only from the pure exponential \(\exp(\beta_{\bar{\lambda}}\varepsilon_n)\), so
\[
E_{r,0}(\bar{\lambda})=\frac{\beta_{\bar{\lambda}}^r}{r!}.
\]
Combining the above expansions gives the stated formula for \(R_{n,\bar{\lambda}}^{t_j}\).
\end{proof}

Applying Lemmas \ref{lemma:dim} and \ref{lemma:ratio-exp}, we now have, for any positive integer $k$,
\begin{alignat}{2}\label{eqn:bigsum}
    \Phi_n &= \frac{n^{tj}e^{-tf_n}}{(tj)!}\,\sum_{\bar{\lambda}\vdash tj} d_{\bar{\lambda}}\chi^{\bar{\lambda}}_{(j^t)}
\left(
\sum_{a=0}^{k}\sum_{r=0}^{k}\sum_{s=0}^{r}
\frac{D_a(\bar{\lambda})E_{r,s}(\bar{\lambda})}{n^a}\,
\varepsilon_n^r\,\frac{1}{f_n^{\,s}}
+o(\varepsilon_n^k)
\right),
\end{alignat}
where $D_0(\bar{\lambda})=1$ independent of $\bar{\lambda}$.

The next lemma is applied to show that the first non-negligible contribution in \eqref{eqn:bigsum} comes from the term where $a,s=0$ and $r=t(j-1)$. Setting $k=t(j-1)$ will then allow us to finish the proof of Theorem \ref{thm:rtransp}.

\begin{lemma}\label{lem:crux}
Let
\[
\beta_{\bar{\lambda}}=\frac{\kappa(\bar{\lambda})-tj(tj+3)}{2j},
\]
where $\kappa(\bar{\lambda})=\sum_{i}\bar{\lambda}_i(\bar{\lambda}_i -2i + 1)$.
Define the linear operator on $\Lambda_{tj}$ by
\[
\mathcal B_{t,j}=\frac1j\Delta-\frac{tj(tj+3)}{2j}\,I,
\]
where $\Delta$ is the cut-and-join operator and $\Delta s_{\bar{\lambda}}=\frac{\kappa(\bar{\lambda})}2 s_{\bar{\lambda}}$. Then
\[
\mathcal B_{t,j}s_{\bar{\lambda}}=\beta_{\bar{\lambda}} s_{\bar{\lambda}}.
\]
Moreover, if
\[
M_q=\sum_{\bar{\lambda}\vdash tj} d_{\bar{\lambda}}\,\chi^{\bar{\lambda}}_{(j^t)}\,\beta_{\bar{\lambda}}^q,
\]
then
\[
M_q=0\quad\text{for }q<t(j-1),
\qquad
M_{t(j-1)}=\frac{(tj)!\,(t(j-1))!}{(j!)^t}.
\]
\end{lemma}

\begin{proof}
Since $\mathcal B_{t,j}$ is diagonal in the Schur basis with eigenvalue $\beta_{\bar{\lambda}}$, the definition of $M_q$ may be rewritten as
\[
M_q=\left\langle \mathcal B_{t,j}^q\,p_{j^t},\,p_{1^{tj}}\right\rangle,
\]
because
\[
p_{j^t}=\sum_{\bar{\lambda}\vdash tj}\chi^{\bar{\lambda}}_{(j^t)}\,s_{\bar{\lambda}},
\qquad
p_{1^{tj}}=\sum_{\bar{\lambda}\vdash tj}d_{\bar{\lambda}}\,s_{\bar{\lambda}},
\]
and the Schur functions are orthonormal for the Hall inner product.

Write the cut-and-join operator as
\[
\Delta=\Delta_+ + \Delta_-,
\]
as in Subsection \ref{subsec: symfuns}.

Recall that the operator $\Delta_+$ increases the number of parts of a power-sum monomial by $1$, while $\Delta_-$ decreases it by $1$. The scalar term in $\mathcal B_{t,j}$ preserves the number of parts.

Now $p_{j^t}$ has $t$ parts and $p_{1^{tj}}$ has $tj=t+t(j-1)$ parts. Hence any term contributing to
\[
\left\langle \mathcal B_{t,j}^q\,p_{j^t},\,p_{1^{tj}}\right\rangle
\]
must increase the number of parts by exactly $t(j-1)$. Since each application of $\mathcal B_{t,j}$ can increase the number of parts by at most $1$, it follows immediately that 
\[
\left\langle \mathcal B_{t,j}^q\,p_{j^t},\,p_{1^{tj}}\right\rangle=0
\] 
for $q<t(j-1)$.

For $q=t(j-1)$, every one of the $t(j-1)$ applications must contribute a split coming from $\Delta_+$; no scalar term can appear, and no $\Delta_-$ term can appear, because either would force the final number of parts to be strictly less than $tj$. It follows that the only contributing terms are those in which every one of the $q$ factors is coming from the $\frac1j\Delta_+$ part of $\mathcal B_{t,j}$. In particular,
\[
[p_{1^{tj}}]\;\mathcal B_{t,j}^{t(j-1)} p_{j^t}
=
j^{-t(j-1)}\,[p_{1^{tj}}]\;\Delta_+^{t(j-1)} p_{j^t}.
\]

Because $\Delta_+$ is a derivation, we have
\[
e^{r\Delta_+}(f g)=e^{r\Delta_+}(f)\,e^{r\Delta_+}(g)
\]
for all symmetric functions $f,g$. Applying this with $f=g=p_j$ and then iterating gives
\[
e^{r\Delta_+}(p_{j^t})=\bigl(e^{r\Delta_+}p_j\bigr)^t.
\]
Taking the coefficient of $r^{t(j-1)}/(t(j-1))!$ yields the multinomial expansion
\[
\Delta_+^{t(j-1)}(p_{j^t})
=
\sum_{q_1+\cdots+q_t=t(j-1)}
\frac{(t(j-1))!}{q_1!\cdots q_t!}
\prod_{i=1}^t \Delta_+^{q_i}p_j.
\]

We now determine which compositions $(q_1,\dots,q_t)$ can contribute to the coefficient of $p_{1^{tj}}$. Since $\Delta_+$ increases the number of parts by exactly $1$ each time, the support of $\Delta_+^a p_j$ consists of monomials with at most $a+1$ parts. To obtain $p_{1^j}$ from $p_j$, one needs exactly $j-1$ applications of $\Delta_+$, so the coefficient of $p_{(1^j)}$ in $\Delta_+^a p_j$ is zero unless $a=j-1$.

Moreover, the coefficient of $p_{1^j}$ in $\Delta_+^{\,j-1}p_j$ is the classical number of ordered minimal transposition factorizations of a fixed $j$-cycle, which by a classical result of Dénes \cite{D59} is
\[
[p_{1^j}]\Delta_+^{\,j-1}p_j=j^{j-2}.
\]
Therefore the only nonzero contribution to $[p_{1^{tj}}]\Delta_+^{t(j-1)} p_{j^t}$ comes from the unique choice
\[
q_1=\cdots=q_t=j-1.
\]
For that choice,
\[
[p_{1^{tj}}]\Delta_+^{t(j-1)} p_{j^t}
=
\frac{(t(j-1))!}{((j-1)!)^t}
\left([p_{1^j}]\Delta_+^{j-1}p_j\right)^t
=
\frac{t(j-1)!}{((j-1)!)^t}\,j^{t(j-2)}.
\]

Since $\langle p_{1^{tj}},p_{1^{tj}}\rangle=(tj)!$, the Hall inner product with $p_{1^{tj}}$ is obtained by multiplying the coefficient of $p_{1^{tj}}$ by $(tj)!$. Hence
\[
\left\langle \Delta_+^{t(j-1)} p_{j^t},\; p_{1^{tj}}\right\rangle
=
(tj)!\,\frac{(t(j-1))!}{((j-1)!)^t}\,j^{t(j-2)}.
\]
Multiplying by the factor $j^{-t(j-1)}$ coming from $\mathcal B_{t,j}$ gives
\[
\left\langle \mathcal B_{t,j}^{t(j-1)} p_{j^t},\; p_{1^{tj}}\right\rangle
=
(tj)!\,\frac{(t(j-1))!}{((j-1)!)^t}\,j^{t(j-2)-t(j-1)}
=
(tj)!\,\frac{(t(j-1))!}{((j-1)!)^t}\,j^{-t}.
\]
Since $j^t(j-1)!^t=(j!)^t$, this simplifies to
\[
\left\langle \mathcal B_{t,j}^{t(j-1)} p_{j^t},\; p_{1^{tj}}\right\rangle
=
\frac{(tj)!(t(j-1))!}{(j!)^t}
\]
\end{proof}

We now have all of the pieces to achieve the necessary cancellations in \eqref{eqn:bigsum} and prove Theorem \ref{thm:rtransp}. 

\begin{proof}[Proof of Theorem \ref{thm:rtransp}]
Begin by setting $k=t(j-1)$ and write, starting from Equation (\ref{eqn:bigsum}):

\begin{alignat}{2}
    \Phi_n &= \frac{n^{tj}e^{-tf_n}}{(tj)!}\,\sum_{\bar{\lambda}\vdash tj} d_{\bar{\lambda}}\chi^{\bar{\lambda}}_{(j^t)}
\left(
\sum_{a=0}^{k}\sum_{r=0}^{k}\sum_{s=0}^{r}
\frac{D_a(\bar{\lambda})E_{r,s}(\bar{\lambda})}{n^a}\,
\varepsilon_n^r\,\frac{1}{f_n^{\,s}}
+o(\varepsilon_n^k)
\right) \nonumber \\
&= \frac{n^{tj}e^{-tf_n}}{(tj)!}\,
\left(
\sum_{a=0}^{k}\sum_{r=0}^{k}\sum_{s=0}^{r}\sum_{\bar{\lambda}\vdash tj} d_{\bar{\lambda}}\chi^{\bar{\lambda}}_{(j^t)}
D_a(\bar{\lambda})E_{r,s}(\bar{\lambda})\frac{\varepsilon_n^r}{f_n^s n^{a}}
+o(\varepsilon_n^k)
\right)\label{eq:startpoint}.
\end{alignat}

Notice that the sum 
\[
\sum_{\bar{\lambda}\vdash tj} d_{\bar{\lambda}}\chi^{\bar{\lambda}}_{(j^t)}
D_a(\bar{\lambda})E_{r,s}(\bar{\lambda}) = \bigl\langle \mathcal D_a\,\mathcal E_{r,s} p_{j^t},\;p_{1^{tj}}\bigr\rangle
\]
is bounded independently of $n$ for any nonnegative integers $a,r,s$ since there are only finitely many partitions of $tj$. Furthermore, since $(j^t)$ has $t$ parts and $(1^{tj})$ has $tj$ parts, we know immediately that
\[
\sum_{\bar{\lambda}\vdash tj} d_{\bar{\lambda}}\chi^{\bar{\lambda}}_{(j^t)}
D_a(\bar{\lambda})E_{r,s}(\bar{\lambda})
= \bigl\langle \mathcal D_a\,\mathcal E_{r,s} p_{j^t},\;p_{1^{tj}}\bigr\rangle=0
\]
whenever $a+r < t(j-1) = k$ by the length filtration conclusions of Lemmas \ref{lemma:dim} and \ref{lemma:ratio-exp}. Hence all terms of total $\varepsilon_n$-degree strictly less than $k$ in \eqref{eq:startpoint} vanish after summation.

It remains to consider the terms with $a+r\ge k$. We have the following cases.\\

\noindent\textit{Case 1: $a+r>k$.} If $a+r>k$, then the contributing term, 
\[
\sum_{\bar{\lambda}\vdash tj} d_{\bar{\lambda}}\chi^{\bar{\lambda}}_{(j^t)}
D_a(\bar{\lambda})E_{r,s}(\bar{\lambda})\frac{\varepsilon_n^r}{f_n^s n^{a}},
\]
is $o(\varepsilon_n^k)$ regardless of $s$. This follows from
\[
\frac{\varepsilon_n^r}{f_n^s n^a}=\frac{f_n^{r-s}}{n^{a+r}}=o\left(\frac{f_n^k}{n^k}\right)\qquad \text{when $a+r>k$}.
\]\\

\noindent\textit{Case 2: $a+r=k$ and $a>0$.} Assume $a+r=k$ and $a>0$. Since 
\[
\mathcal D_a(\mathsf{P}_{\leq \ell}) \subseteq \mathsf{P}_{\leq \ell+a-1} \qquad \text{when $a>0$,}
\]
by Lemmas \ref{lemma:dim} and \ref{lemma:ratio-exp} we have that
\[
\mathcal D_a \mathcal E_{r,s}p_{j^t} \subseteq \mathsf{P}_{\leq t+r+a-1} = \mathsf{P}_{\leq t+k-1} = \mathsf{P}_{\leq tj-1}.
\] 
But clearly $p_{1^{tj}}\notin \mathsf{P}_{\leq tj-1}$, so we conclude that 
\[
\sum_{\bar{\lambda}\vdash tj} d_{\bar{\lambda}}\chi^{\bar{\lambda}}_{(j^t)}
D_a(\bar{\lambda})E_{r,s}(\bar{\lambda})=\bigl\langle \mathcal D_a\,\mathcal E_{r,s} p_{j^t},\;p_{1^{tj}}\bigr\rangle=0.
\]\\

\noindent\textit{Case 3: $r=k$, $a=0$, $s=0$.} In the case that $a=0$, $r=k$, and $s=0$ since $D_0(\bar{\lambda})=1$ and $E_{k,0}(\bar{\lambda}) = \frac{\beta_{\bar{\lambda}}^{k}}{k!}$, we have
\[
\sum_{\bar{\lambda}\vdash tj} d_{\bar{\lambda}}\chi^{\bar{\lambda}}_{(j^t)}
D_0(\bar{\lambda})E_{k,0}(\bar{\lambda})=\frac{M_{k}}{k!}=\frac{(tj)!}{(j!)^t}
\]
by Lemma \ref{lem:crux}.\\

\noindent\textit{Case 4: $a+r=k$ and $s>0$.} If $a+r=k$ and $s\ge 1$, then the corresponding term has an additional factor $f_n^{-s}$. We have
\[
\frac{\varepsilon_n^r}{f_n^sn^a} = \frac{f_n^{r-s}}{n^k} = o\left(\frac{f_n^k}{n^k}\right) \qquad \text{since $r-s=k-a-s<k$},
\]
and hence the term is $o(\varepsilon_n^k)$.\\

In summary, we obtain
\begin{alignat*}{2}
\Phi_n
&=
\frac{n^{tj} e^{-tf_n}}{(tj)!}
\left(
\frac{(tj)!}{(j!)^t}\varepsilon_n^{t(j-1)}
+o\!\left(\varepsilon_n^{t(j-1)}\right)
\right)\\
&= \frac{n^{tj}e^{-tf_n}\varepsilon_n^{t(j-1)}}{(j!)^t}(1+o(1)).
\end{alignat*}

Now
\[
n^{tj} e^{-tf_n}\varepsilon_n^{t(j-1)}
=
n^{tj}\,e^{-t(\log n+(j-1)\log\log n+c)}
\left(\frac{\log n+(j-1)\log\log n+c}{n}\right)^{t(j-1)},
\]
which simplifies to
\[
e^{-tc}\,(1+o(1)).
\]
Hence
\[
\lim_{n\to\infty}\Phi_n=\frac{e^{-tc}}{(j!)^t}.
\]
\end{proof}

\section{Bratteli diagrams and proof of Theorem \ref{thm:A}}\label{bratteli}

The main goal of this section is to prove Theorem \ref{thm:A}. Let us first briefly recall the motivation. We want to compute the moments of the class function 
 \[
 a_j(\sigma) = \{\text{\# of $j$-cycles of $\sigma$}\}
 \]
 over probability measures of interest. This in turn will allow us to prove Poisson limits for finite $j$-cycle count distributions of permutations obtained from random walks on $S_n$.

Translating from characters to representations, we seek to decompose tensor powers of the virtual $S_n$-representation $\rho_{a_j}$ having character $a_j$ as a formal sum of irreducible representations. Note that the $j=1$ case (fixed points) corresponds to decomposing tensor powers of the $n$-dimensional defining representation (identical to $\rho_{a_1}$) as a direct sum of irreducible representations, which is achievable by several methods \cite{F24, D16, BHH16}. Ding \cite{D16} uses a generating function approach. Fulman \cite{F24} uses connections with card shuffling. Benkart, Halverson, and Harmon \cite{BHH16} provide the approach that inspire our formula and its proof. They introduce restriction-induction Bratteli diagrams that capture the notion of ``removing and adding a box" to Young diagrams in all possible ways at increasing levels of a Bratteli diagram, which in turn corresponds to taking successive tensor products with the defining representation. 

Our method for computing the decomposition of the $r$th tensor power of $\rho_{a_j}$ for fixed $j$ between 1 and $\frac{n}{2r}$ uses the main idea that, roughly, tensoring with the closely related virtual representation $\rho_{\psi_j}$ corresponds to ``removing and adding rim-$j$ hooks" to Young diagrams in all possible ways. The MN$_j$ restriction and induction rules we introduce in Definition \ref{def:mnj} can be seen as an alternate way (allowing signs in accordance with the Murnaghan-Nakayama rule) to generalize restriction and induction of the group-subgroup pair $(S_n, S_{n-j})$ in the ring of virtual $S_n$ representations from the well-known rule \ref{eq:RI} for $(S_n,S_{n-1})$. 

Subsection \ref{subsec:BratSn} provides useful background on restriction-induction Bratteli diagrams for the pair $(S_n,S_{n-1})$, which Subsection \ref{subsec:mnj} extends to our new tool, MN$_j$ restriction-induction Bratteli diagrams for the pair $(S_n, S_{n-j})$. Finally, Subsection \ref{proofsec} completes proof of Theorem \ref{thm:A}.

\subsection{Restriction-induction Bratteli diagrams}\label{subsec:BratSn}\hfill\\

In this subsection, we explain the machinery behind the restriction-induction Bratteli diagrams that can be used to compute tensor powers of the defining representation of $S_n$.

We immediately specialize to the case of the group-subgroup pair
$(S_n, S_{n-1})$ for our purposes. Further details including a discussion of restriction-induction Bratteli diagrams for general group-subgroup pairs can be found in \cite{BH17} and \cite{BHH16}.   

The \emph{restriction-induction Bratteli diagram} for the pair $(S_n, S_{n-1})$  is an infinite rooted tree $\mathcal{B}(S_n, S_{n-1})$ whose vertices are organized into rows labeled by half integers $\ell$  in $\half\N$.  For $\ell =r \in \N$, the vertices on row
$r$ are partitions of $n$, and the vertices on row $\ell = r+\half$ are partitions of $n-1$. 

 The edges of $\mathcal{B}(S_n, S_{n-1})$ are constructed from the restriction and induction rules for $(S_n, S_{n-1})$, which are well known (see Sagan \cite{S01} for example):
\begin{equation}\label{eq:RI}  \mathsf{Res}^{S_n}_{S_{n-1}}(S^\lambda)  = \bigoplus_{\mu = \lambda-\square} S^\mu, 
\qquad  \mathsf{Ind}^{S_{n+1}}_{S_{n}}(S^\lambda)  = \bigoplus_{\nu = \lambda+\square} S^\nu, \end{equation}
where the first sum is over all partitions $\mu$ of $n-1$ obtained from $\lambda$ by removing a box from an inner corner of the diagram of
$\lambda$, and
the second sum is over all partitions $\nu$ of $n+1$ obtained by adding a box to an outer corner of $\lambda$.  

Applying these rules to the trivial $S_n$-module, we see that
\[ \Ind^{S_n}_{S_{n-1}}\big(\Res^{S_n}_{S_{n-1}}(S^{(n)})\big)  = 
 \Ind^{S_n}_{S_{n-1}}\big(S^{(n-1)}\big) =
 S^{(n)} \oplus S^{(n-1,1)} \cong \varrho_n.\] 

Crucially, the restriction and induction operators on $S_n$ modules satisfy the property that for any (virtual) $S_n$-module $X$, 
\begin{equation}\label{eqn:prop1} 
\Ind^{S_n}_{S_{n-1}}\big(\Res^{S_n}_{S_{n-1}}(X)\big) \cong X \otimes \varrho_n.
\end{equation}

 We track the number of paths from $(n)$ in level 0 to $\lambda$ at integer level $\ell = r$ in $\mathcal{B}(S_n, S_{n-1})$ by labeling vertex $(\lambda, r)$ with coefficient $c_{\lambda,r}$. The $c_{\lambda,r}$ can be computed recursively using 
 
 \[c_{\lambda,r} = \sum_{\mu: \lambda/\mu = \square} c_{\mu, r-\half} \qquad \text{and} \qquad
c_{\mu,r+\half} = \sum_{\lambda: \lambda/\mu = \square} c_{\lambda, r}.\]

An example of a restriction-induction Bratteli diagram for the pair $(S_7,S_6)$ is illustrated in Figure \ref{fig:Sbratteli}. See \cite{BH17, BHH16} for more examples.

 By construction, we have the following theorem from \cite{BHH16}, which illustrates the power of the restriction-induction Bratteli diagram as a tool for computing tensor powers of the defining representation of $S_n$.
 
\begin{theorem}
    For $r,n \in \N$ and $n\geq 1$, for all $\lambda \vdash n$, the signed multiplicity of irrep $S^\lambda$ in $\varrho_n^{\otimes r}$ is given by
\begin{equation*}
m_{\lambda,r}^1  =  
\big \vert \big\{\text{paths in $\mathcal{B}(S_n,S_{n-1})$ from $(n)$ in level 0 to $\lambda$ in level $r$}\big\} \big \vert.
\end{equation*}
\end{theorem}

\Yboxdim{8pt}  
\Ylinethick{.6pt}  
\begin{figure}[H]
\[
\begin{tikzpicture}[line width=.5pt,xscale=0.15,yscale=0.25]  
\path (-10,0)  node[anchor=west]  {$\boldsymbol{\ell=0}$};
\path (0,0)  node[anchor=west] (S7-0) {\yng(7)};
\draw (S7-0) node[below=4pt, blue] {\nmm 1};
\path (-10,-5)  node[anchor=west]  {$\boldsymbol{\ell=\frac{1}{2}}$};
\path (10,-5)  node[anchor=west] (S6-1) {\yng(6)};
\draw (S6-1) node[below=4pt,blue] {\nmm 1};
\path (-10,-10)  node[anchor=west]  {$\boldsymbol{\ell=1}$};
\path (0,-10)  node[anchor=west] (S7-2) {\yng(7)};
\path (19,-10)  node[anchor=west] (S61-2) {\yng(6,1)};
\draw (S7-2) node[below=4pt,blue] {\nmm 1};
\draw (S61-2) node[below=6pt, right=-2pt,blue] {\nmm 1};
\path (-10,-15)  node[anchor=west]  {$\boldsymbol{\ell=\frac{3}{2}}$};
\path (10,-15)  node[anchor=west] (S6-3) {\yng(6)};
\path (29,-15)  node[anchor=west] (S51-3) {\yng(5,1)};
\draw (S6-3) node[below=4pt, blue] {\nmm 2};
\draw (S51-3) node[below=6pt,right=2pt,blue] {\nmm 1};
\path (-10,-20)  node[anchor=west]  {$\boldsymbol{\ell=2}$};
\path (0,-20)  node[anchor=west] (S7-4) {\yng(7)};
\path (19,-20)  node[anchor=west] (S61-4) {\yng(6,1)};
\path (34,-20)  node[anchor=west] (S52-4) {\yng(5,2)};
\path (49,-20)  node[anchor=west] (S511-4) {\yng(5,1,1)};
\draw (S7-4) node[below=4pt, blue] {\nmm 2};
\draw (S61-4) node[below=6pt, right=-2pt,blue] {\nmm 3};
\draw (S52-4) node[below=4pt, right=8pt, blue] {\nmm 1};
\draw (S511-4) node[below=4pt, right=2pt, blue] {\nmm 1};

\path  (S7-0) edge[black,thick] (S6-1);
\path  (S7-2) edge[black,thick] (S6-1);
\path  (S61-2) edge[black,thick] (S6-1);
\path  (S7-2) edge[black,thick] (S6-3);
\path  (S61-2) edge[black,thick] (S6-3);
\path  (S61-2) edge[black,thick] (S51-3);
\path  (S7-4) edge[black,thick] (S6-3);
\path  (S61-4) edge[black,thick] (S6-3);
\path  (S61-4) edge[black,thick] (S51-3);
\path  (S52-4) edge[black,thick] (S51-3);
\path  (S511-4) edge[black,thick] (S51-3);
\end{tikzpicture}
\]
    \caption{Levels  \ $\ell = 0,\half,1,\frac{3}{2},2$ \  of the restriction-induction Bratteli diagram for the pair $(S_7, S_6)$. Reading off coefficients and indexing partitions from the $\ell=2$ level of the diagram gives the irreducible decomposition of $\varrho_7^{\otimes 2}$. \label{fig:Sbratteli}}
\end{figure}

\subsection{MN$_j$ restriction-induction Bratteli diagrams}\label{subsec:mnj}\hfill\\

With the goal of taking advantage of a property like (\ref{eqn:prop1}) to compute tensor powers of $\rho_{\psi_j}$ for general fixed $j$ by counting paths through a Bratteli diagram, we introduce Murnaghan-Nakayama-$j$, ($\textrm{MN}_j$) restriction and induction, defined on the Specht modules as follows.

\begin{definition}\label{def:mnj}
Let $n,j \in \ZZ_{\geq 1}$ with $n>j$. Using notation of Sagan \cite{S01}, we write $ll(\lambda/\mu)$ to denote the leg length (the number of rows minus 1) of rim hook $\lambda/\mu$. For $\lambda \vdash n$,
\[R_j^{\textrm{MN}}(S^\lambda) = \sum_{\substack{\mu\vdash n-j \\ \lambda/\mu \textrm{ a rim-$j$ hook}}} (-1)^{ll(\lambda/\mu)}S^\mu.\]

\noindent
For $\mu \vdash n-j$,
\[I_j^{\textrm{MN}}(S^\mu) = \sum_{\substack{\lambda\vdash n \\ \lambda/\mu \textrm{ a rim-$j$ hook}}} (-1)^{ll(\lambda/\mu)}S^\lambda.\]
\end{definition}

These operations extend linearly in the complex representation ring $R(S_n)\otimes_{\mathbb{Z}}\mathbb{C}$ and coincide with ordinary restriction and induction when $j=1$. Corresponding statements in terms of symmetric functions amount to the following two dual restatements of the Murnaghan-Nakayama rule:
\[
p_j^\perp s_\lambda = \sum_{\substack{\mu\vdash n-j \\ \lambda/\mu \textrm{ a rim-$j$ hook}}} (-1)^{ll(\lambda/\mu)}s_\mu
\]
and
\[
p_j s_\mu = \sum_{\substack{\lambda\vdash n \\ \lambda/\mu \textrm{ a rim-$j$ hook}}} (-1)^{ll(\lambda/\mu)}s_\lambda,
\]
where $p_j^\perp$ is the adjoint of multiplication by power sum $p_j$ with respect to the Hall inner product.

To use the Bratteli diagram machinery developed in Subsection \ref{subsec:BratSn}, we want to say that taking a tensor product of (virtual) $S_n$-representation $X$ with $\rho_{\psi_j}$ is equivalent to performing $\textrm{MN}_j$-restriction followed by $\textrm{MN}_j$-induction on $X$. Towards this claim, we verify the following ``tensor identity''.

\begin{theorem}\label{thm:key}
For $X$ a virtual $S_n$-module and $Y$ a virtual $S_{n-j}$-module, we have an isomorphism of virtual $S_n$-modules
\begin{equation}\label{eqn:tensor}
\mathrm{I}_j^{\mathrm{MN}}\!\big(\mathrm{R}_j^{\mathrm{MN}}(X)\otimes Y\big)
\;\cong\;
X\otimes \mathrm{I}_j^{\mathrm{MN}}(Y).
\end{equation}
\end{theorem}

\begin{proof}
Translating to symmetric functions via the Frobenius characteristic map, the identity we seek to prove is
\begin{equation}\label{eqn:sym}
    p_j \big(p_j^\perp f * g\big)
=
f * (p_j g),
\end{equation}
for symmetric functions $f$ and $g$ of degrees $n$ and $n-j$ respectively with $n>j$. The operators $p_j$, $p_j^\perp$ are linear and the Kronecker product $\ast$ is bilinear. So it suffices to consider $f = p_\alpha$ and $g = p_\beta$ with $\alpha \vdash n$, $\beta \vdash n-j$ since the power sum symmetric functions form a $\mathbb{Q}$-basis for the ring of symmetric functions. Recall that for partitions $\mu$ and $\nu$, the Kronecker product satisfies
\[
p_\mu \ast p_\nu = \delta_{\mu,\nu} \, z_\mu \, p_\mu,
\]
where $z_\mu = \prod_i i^{m_i(\mu)} m_i(\mu)!$ and $m_i(\mu)$ is the signed multiplicity of the part $i$ in $\mu$. Notice that if $\alpha$ does not contain $j$ as a component, then both sides of (\ref{eqn:sym}) are 0. So assume $\alpha = j \cup \gamma$ where $\cup$ denotes concatenation of parts of partitions as multisets (e.g. $(4,2,1) \cup (3,3,1) = (4,3,3,2,1,1))$. In this case, we have the following claim.\\

\noindent\textbf{Claim: }
Assuming $\alpha = j \cup \gamma$,
    \[p_j^\perp p_\alpha = \frac{z_\alpha}{z_\gamma}p_\gamma.\]

\noindent\textit{Proof of claim.} Let $r = m_j(\alpha)$ be the signed multiplicity of the part $j$ in $\alpha$. 
Let $\gamma$ be the partition obtained by removing exactly one copy of $j$ from $\alpha$, so $|\gamma| = |\alpha| - j$ and $\gamma$ contains $r-1$ copies of $j$.  

By definition of $p_j^\perp$ as the adjoint of multiplication by $p_j$ with respect to the Hall inner product (see Macdonald \cite[p. 76]{M95}):
\[
\langle p_j^\perp p_\alpha, p_\gamma \rangle = \langle p_\alpha, p_j p_\gamma \rangle = \langle p_\alpha, p_\alpha \rangle = z_\alpha.
\]

On the other hand, if $p_j^\perp p_\alpha = c \, p_\gamma$, then $\langle p_j^\perp p_\alpha, p_\gamma \rangle = c \, \langle p_\gamma, p_\gamma \rangle = c \, z_\gamma$ since the Hall inner product is diagonal in the power sum basis. Therefore \[c = \frac{z_\alpha}{z_\gamma},\] and the claim follows. \qed\\

Returning to the proof of (\ref{eqn:sym}), let $f = p_\alpha$ and $g = p_\beta$ and further assume that $\alpha = j \cup \gamma$ ($\alpha$ contains $j$ as a component). We have, for the left hand side of (\ref{eqn:sym}),

\[
p_j \big( (p_j^\perp p_\alpha) \ast p_\beta \big) 
= p_j \Big( \frac{z_\alpha}{z_\gamma} p_\gamma \ast p_\beta \Big).
\]

By the Kronecker product formula, $p_\gamma \ast p_\beta = 0$ unless $\beta = \gamma$, in which case $p_\gamma \ast p_\gamma = z_\gamma \, p_\gamma$. Therefore,
\[
p_j \big( (p_j^\perp p_\alpha) \ast p_\beta \big) =
\begin{cases}
z_\alpha \, p_j p_\gamma = z_\alpha \, p_\alpha, & \text{if } \beta = \gamma, \\
0, & \text{otherwise.}
\end{cases}
\]

Meanwhile, on the right hand side of (\ref{eqn:sym}), we have
\[
p_\alpha \ast (p_j p_\beta) = p_\alpha \ast p_{j \cup \beta}.
\]
If $\beta = \gamma$, then
\[
p_\alpha \ast (p_j p_\beta) = p_\alpha \ast p_\alpha = z_\alpha \, p_\alpha.
\]
If $\beta \neq \gamma$ then $j\cup \beta \neq j \cup \gamma = \alpha$ and the right hand side is 0. The two sides of (\ref{eqn:sym}) thus match exactly.
\end{proof}

With Theorem \ref{thm:key} in hand, we set $Y$ equal to the trivial representation of $S_{n-j}$ in (\ref{eqn:tensor}) and obtain the key identity: for any virtual $S_n$-module $X$, we have 
\begin{equation}\label{eqn:key}
    I_j^{\mathrm{MN}}(R_j^{\mathrm{MN}}(X)) \cong X \otimes \rho_{\psi_j},
\end{equation}
making use of the fact that 
\[
\rho_{\psi_j} \cong I_j^{\mathrm{MN}}(S^{(n-j)}).
\]

Letting $X$ be the trivial $S_n$-module at the level 0 root, we can now compute tensor powers of $\rho_{\psi_j}$ using a $\mathrm{MN}_j$ restriction-induction Bratteli diagram for the pair $(S_n, S_{n-j})$, which we denote by $\mathcal{MN}_j(S_n, S_{n-j})$. In the example diagram given by Figure \ref{fig:mbratteli1} below, we compute a decomposition of $\rho_{\psi_2}^{\otimes 2}$ valid for $n\geq 8$.

\Yboxdim{6pt}  
\Ylinethick{.6pt} 
\begin{figure}[H]
\[
\begin{tikzpicture}[line width=.5pt,xscale=0.12,yscale=0.2]  
\path (-15,0)  node[anchor=west]  {\boldsymbol{$\ell = 0$}};
\path (10,0)  node[anchor=west] (S8-0) {\yng(8)};
\draw (S8-0) node[below=4pt, black] {\color{blue}\nmm 1};
\path (-15,-7.5)  node[anchor=west]  {\boldsymbol{$\ell = \half$}};
\path (25,-7.5)  node[anchor=west] (S6-1) {\yng(6)};
\draw (S6-1) node[below=4pt,black] {\color{blue}\nmm 1};
\path (-15,-15)  node[anchor=west]  {\boldsymbol{$\ell = 1$}};
\path (10,-15)  node[anchor=west] (S8-2) {\yng(8)};
\path (30,-15)  node[anchor=west] (S62-2) {\yng(6,2)};
\path (50,-15)  node[anchor=west] (S611-2) {\yng(6,1,1)};
\draw (S8-2) node[below=8pt, left=1pt, black] {\color{blue}\nmm 1};
\draw (S62-2) node[below=3pt,black] {\color{blue}\nmm 1};
\draw (S611-2) node[below=1pt,black] {\color{blue}\nmm -1};
\path (-15,-22.5)  node[anchor=west]  {\boldsymbol{$\ell = \frac{3}{2}$}};
\path (15,-22.5)  node[anchor=west] (S6-3) {\yng(6)};
\path (35,-22.5)  node[anchor=west] (S42-3) {\yng(4,2)};
\path (55,-22.5)  node[anchor=west] (S411-3) {\yng(4,1,1)};
\draw (S6-3) node[below=8pt, left=13pt, black] {\color{blue}\nmm 3};
\draw (S42-3) node[below=11pt, left=1pt, black] {\color{blue}\nmm 1};
\draw (S411-3) node[below=2pt,black] {\color{blue}\nmm -1};
\path (-15,-30)node[anchor=west]  {\boldsymbol{$\ell = 2$}};
\path (-3,-30)  node[anchor=west] (S8-4) {\yng(8)};
\path (12,-30)  node[anchor=west] (S62-4) {\yng(6,2)};
\path (25,-30)  node[anchor=west] (S611-4) {\yng(6,1,1)};
\path (38,-30)  node[anchor=west] (S44-4) {\yng(4,4)};
\path (65,-30)  node[anchor=west] (S431-4) {\yng(4,3,1)};
\path (56,-30)  node[anchor=west] (S422-4) {\yng(4,2,2)};
\path (47,-30)  node[anchor=west] (S4211-4) {\yng(4,2,1,1)};
\path (74,-30)  node[anchor=west] (S41111-4) {\yng(4,1,1,1,1)};
\draw (S8-4) node[below=5pt, black] {\color{blue}\nmm 3};
\draw (S62-4) node[below=5pt, black] {\color{blue}\nmm 4};
\draw (S611-4) node[below=5pt, black] {\color{blue}\nmm -4};
\draw (S44-4) node[below=5pt, black] {\quad\color{blue}\nmm 1};
\draw (S431-4) node[below=5pt, black] {\quad\color{blue}\nmm -1};
\draw (S422-4) node[below=5pt, black] {\quad\color{blue}\nmm 2};
\draw (S4211-4) node[below=5pt, black] {\quad\color{blue}\nmm -1};
\draw (S41111-4) node[below=5pt, black] {\quad\color{blue}\nmm 1};

\path  (S8-0) edge[black,thick] (S6-1);
\path  (S8-2) edge[black,thick] (S6-1);
\path  (S62-2) edge[black,thick] (S6-1);
\path  (S611-2) edge[red,thick] (S6-1);
\path  (S8-2) edge[black,thick] (S6-3);
\path  (S62-2) edge[black,thick] (S6-3);
\path  (S62-2) edge[black,thick] (S42-3);
\path  (S611-2) edge[red,thick] (S6-3);
\path  (S611-2) edge[black,thick] (S411-3);
\path  (S8-4) edge[black,thick] (S6-3);
\path  (S62-4) edge[black,thick] (S6-3);
\path  (S62-4) edge[black,thick] (S42-3);
\path  (S44-4) edge[black,thick] (S42-3);
\path  (S611-4) edge[red,thick] (S6-3);
\path  (S611-4) edge[black,thick] (S411-3);
\path  (S4211-4) edge[red,thick] (S42-3);
\path  (S422-4) edge[black,thick] (S42-3);
\path  (S422-4) edge[red,thick] (S411-3);
\path  (S431-4) edge[black,thick] (S411-3);
\path  (S41111-4) edge[red,thick] (S411-3);

\end{tikzpicture}
\]
\caption{Levels  \ $0,\half,1,\frac{3}{2},2$ \  of the $\mathrm{MN}_2$ restriction-induction Bratteli diagram for the pair $(S_{8}, S_6)$. Edges drawn in red indicate addition/removal of rim-hooks with odd leg length.\label{fig:mbratteli1}}
\end{figure}

Note the interaction of signed edges with the labeled coefficients in the diagram. When a red edge connects a vertex $\mu$ in level $\ell$ to one, say $\nu$, in level $\ell + \half$, the coefficient of vertex $\mu$ is subtracted from the total in the calculation of the coefficient of $\nu$ instead of added. Concretely, the vertex labeled by $\lambda$ at level $r$ of the diagram is assigned coefficient $m_{\lambda,r}^j$, which can be computed recursively by the rules
\[
m_{\lambda,r}^j = \sum_{\substack{\mu\vdash n-j \\ \lambda/\mu \textrm{ a rim-$j$ hook}}}(-1)^{ll(\lambda/\mu)}m_{\mu,r-\half}^j
\]
and
\[
m^j_{\mu,r+\half} = \sum_{\substack{\lambda\vdash n \\ \lambda/\mu \textrm{ a rim-$j$ hook}}}(-1)^{ll(\lambda/\mu)} m^j_{\lambda, r}.
\]
The stabilization of coefficients in the diagram as $n$ grows is obvious. In the example of Figure \ref{fig:mbratteli1}, the signed multiplicity of $S^{(n-tj,\bar{\lambda})}$ in $\rho_{\psi_j}^{\otimes r}$ (with $0\leq t \leq r$ and $\bar{\lambda}\vdash tj$) is identical for all $n\geq 2rj$. We conclude that for $n\geq 8$, \[\rho_{\psi_2}^{\otimes 2} = 3S^{(n)} + 4S^{(n-2,2)} - 4S^{(n-2,1^2)}+\sum_{i=0}^{3}(-1)^{\min(i,3-i)} S^{(n-4,4-i,1^i)}+2S^{(n-4,2^2)}.\] 

\noindent Equivalently, we have the irreducible character decompositions
\[(\psi_2)^2 = 3 + 4\chi^{(n-2,2)} - 4\chi^{(n-2,1^2)}+\sum_{i=0}^{3}(-1)^{\min(i,3-i)} \chi^{(n-4,4-i,1^i)}+2\chi^{(n-4,2^2)},\]
and 
\[(a_2)^2 = \frac{1}{4}\Big(3 + 4\chi^{(n-2,2)} - 4\chi^{(n-2,1^2)}+\sum_{i=0}^{3}(-1)^{\min(i,3-i)} \chi^{(n-4,4-i,1^i)}+2\chi^{(n-4,2^2)}\Big).\]

By construction of $\mathcal{MN}_j(S_n,S_{n-j})$, we have Theorem \ref{thm:pathcount}, which reduces the proof of Theorem \ref{thm:A} to a path counting problem. The sign of $m_{\lambda,r}^j$ should be negative when an odd number of rim-hooks with odd leg length are added/removed, and positive when this number is even. Note that the number of odd-signed edges in a path in $\mathcal{MN}_j(S_n,S_{n-j})$ from $(n)$ at level 0 to $\lambda$ at level $r$ does not depend on the choice of path.

\begin{theorem}\label{thm:pathcount}
    For $j,r,n \in \N$ with $j\geq 1$, $n\geq 2j$ and any $\lambda \vdash n$, let 
    \[
    S = \big\{\text{paths in $\mathcal{MN}_j(S_n,S_{n-j})$ from $(n)$ in level 0 to $\lambda$ in level $r$}\big\}
    \]
    and $s$ be 0 if the number of odd-signed edges in any path in $S$ is even, else 1. The signed multiplicity of irrep $S^\lambda$ in $\rho_{\psi_j}^{\otimes r}$ is given by
\begin{equation*}
m_{\lambda,r}^j  =  (-1)^{s} \big \vert S \big \vert.
\end{equation*}
\end{theorem} 

\subsection{Proof of Theorem \ref{thm:A}}\label{proofsec}\hfill\\

The work of the previous subsection shows how to recursively compute the (integer-valued) signed multiplicity of irrep $S^\lambda$ in $\rho_{\psi_j}^{\otimes r}$ by counting paths through a diagram. While Theorem \ref{thm:pathcount} makes only the assumption $n\geq 2j$ for which $\rho_{\psi_j}$ is defined on the size of $n$, we will assume throughout this section that $n\geq 2rj$ with the goal of proving Theorem \ref{thm:A}.

\begin{remark} It is quite possible that a compact expression for $m_{\lambda,r}^j$ can be obtained without requiring $n\geq 2rj$; with the goal of this paper being to describe the application of the irreducible decomposition of $\rho_{\psi_j}^{\otimes r}$ to the limiting distributions of cycle statistics of permutations obtained from random walks on $S_n$ as $n$ grows large, we are content to leave it for another day.
\end{remark}

To prove Theorem \ref{thm:A}, we need to show that for $j\geq 1$ and $n\geq 2rj$, the signed multiplicity $m_{\lambda,r}^j$ in $\rho_{\psi_j}^{\otimes r}$ is given by
\begin{equation}\label{eqn:mform}
m^{j}_{\lambda,r}
= \begin{cases}
    \chi^{\bar{\lambda}}_{(j^t)}c_{t,r}^j &\quad \text{when $\lambda = (n-tj,\bar{\lambda})$},\\
    0 &\quad \text{otherwise},
\end{cases}
\end{equation}
where $c_{t,r}^j = \sum_{a=t}^r S(r,a)\binom{a}{t}j^{r-a}$.

Let us analyze the two multiplicative factors of $m_{\lambda,r}^j$ separately at first since it is not clear from Theorem \ref{thm:pathcount} why \eqref{eqn:mform} should hold.

We can compute the magnitude of $\chi^{\bar{\lambda}}_{(j^t)}$ as the number of paths from $(n)$ in level 0 to $(n-tj,\bar{\lambda})$ in level $t$ of $\mathcal{MN}_j(S_n,S_{n-j})$ (noting that $c_{t,t}^j = 1$). This can be equivalently counted as the number of ways to assemble $\bar{\lambda}$ using consecutive rim-$j$ hook additions, which has been previously analyzed in \cite{J78, JK81, S99} and the references therein. Furthermore, the Murnaghan-Nakayama rule guarantees that $\chi^{\bar{\lambda}}_{(j^t)}$ carries the same sign as $m_{\lambda,r}^j$ when $\lambda = (n-tj, \bar{\lambda})$ and $t\leq r$. This sign tracks the number (even or odd) of components with odd leg length in any rim-$j$ filling of $\bar{\lambda}$, or equivalently the number of odd-signed edges in any path from $(n)$ in level 0 to $(n-tj,\bar{\lambda})$ in level $r$ of $\mathcal{MN}_j(S_n,S_{n-j})$. 

Meanwhile, the $c_{t,r}^j$ can be computed by counting paths through $\mathcal{MN}_j(S_n,S_{n-j})$ arriving at a partition $\lambda=(n-tj,\bar{\lambda})$ in level $r$ having the property that $\chi^{\bar{\lambda}}_{(j^t)} = 1$. A sensible choice then is to analyze paths reaching $\lambda = (n-tj,tj)$ in level $r$ from $(n)$ in level 0. Noting that there are no paths through $\mathcal{MN}_j(S_n,S_{n-j})$ started from $(n)$ in level 0 to partitions (in any level) with first row lengths other than $n-tj$ for nonnegative integers $t$, we count paths from $(n)$ in level 0 to $(n-tj,tj)$ in level $r \geq t$ to obtain $m_{\lambda,r}^j$ for $\lambda = (n-tj,tj)$ in a direct combinatorial manner. The following theorem provides a closed formula involving the Stirling numbers of the second kind for this path count.

\begin{theorem}\label{thm:paths}
    Let $j,t,r,n \in \N$, $j \geq 1$, $n\ge 2rj$ and $0\le t\le r$. The number of paths from $(n)$ in level 0 to $(n-tj,tj)$ in level $r$ in the MN$_j$ restriction-induction Bratteli diagram for the pair $(S_n,S_{n-j})$ is given by
    \[c_{t,r}^j = \sum_{a=t}^r S(r,a)\binom{a}{t}j^{r-a},\]
    where $S(r,a)$ denotes the Stirling number of the second kind.
\end{theorem}

\begin{proof}
    Let $p_{t,r}^j$ be the number of paths in $r$ rim-$j$ hook removal-plus-addition steps from $(n)$ to $(n-tj,tj)$ and $c_{t,r}^j$ be the sequence given by Theorem \ref{thm:paths}. We will show that $p_{t,r}^j$ satisfies the same recurrence relation and boundary conditions as $c_{t,r}^j$ and thus counts the same thing. Let us begin with the $c_{t,r}^j$. Using the Stirling recurrence
\(S(r+1,i)=iS(r,i)+S(r,i-1)\), we have:
\[
\begin{aligned}
c_{t,r+1}^j
&=\sum_{i=t}^{r+1} S(r+1,i)\binom{i}{t}j^{r+1-i} \\
&=\sum_{i=t}^{r+1}\big(iS(r,i)+S(r,i-1)\big)\binom{i}{t}j^{r+1-i} \\
&=\sum_{i=t}^{r+1} iS(r,i)\binom{i}{t}j^{r+1-i}
\;+\;\sum_{i=t}^{r+1} S(r,i-1)\binom{i}{t}j^{r+1-i}.
\end{aligned}
\]

\noindent For the second sum, change the index by setting \(k=i-1\):
\begin{equation}\label{eqn:stir1}
\sum_{i=t}^{r+1} S(r,i-1)\binom{i}{t}j^{r+1-i}
=\sum_{k=t-1}^{r} S(r,k)\binom{k+1}{t}j^{r-k}.
\end{equation}
Now use the binomial recurrence \(\binom{k+1}{t}=\binom{k}{t}+\binom{k}{t-1}\) to split (\ref{eqn:stir1}) into
\begin{equation}\label{eqn:stir2}
\sum_{k=t}^{r} S(r,k)\binom{k}{t}j^{r-k}
\;+\;\sum_{k=t-1}^{r} S(r,k)\binom{k}{t-1}j^{r-k}
= c_{t,r}^j + c_{t-1,r}^j.
\end{equation}

\noindent Next treating the first sum in \ref{eqn:stir1}, we have
\[
\sum_{i=t}^{r+1} iS(r,i)\binom{i}{t}j^{r+1-i}
= j\sum_{i=t}^{r} iS(r,i)\binom{i}{t}j^{r-i},
\]
because \(S(r,r+1)=0\). Applying the binomial identity
\[
i\binom{i}{t}=t\binom{i}{t}+(t+1)\binom{i}{t+1},
\]
we obtain
\[
\begin{aligned}
j\sum_{i=t}^{r} iS(r,i)\binom{i}{t}j^{r-i}
&= j\sum_{i=t}^{r}\Big(t\binom{i}{t}+(t+1)\binom{i}{t+1}\Big)S(r,i)j^{r-i} \\
&= j\Big(t\sum_{i=t}^{r} S(r,i)\binom{i}{t}j^{r-i}
+(t+1)\sum_{i=t+1}^{r} S(r,i)\binom{i}{t+1}j^{r-i}\Big) \\
&= j\Big(t\,c_{t,r}^j+(t+1)\,c_{t+1,r}^j\Big).
\end{aligned}
\]

Combining the two pieces (the contribution from the \(iS(r,i)\) part
and the shifted \(S(r,i-1)\) part) yields
\[
c_{t,r+1}^j
= j\big(t\,c_{t,r}^j+(t+1)\,c_{t+1,r}^j\big) + \big(c_{t,r}^j+c_{t-1,r}^j\big).
\]
Collecting terms and shifting gives the recurrence:
\begin{equation}\label{*}
   c_{t,r}^j = c_{t-1,r-1}^j + (j\,t+1)\,c_{t,r-1}^j + j\,(t+1)\,c_{t+1,r-1}^j. 
\end{equation}
We additionally have the boundary conditions:
\begin{itemize}
    \item $c_{r,r}^j = 1$ for $r\geq 0$
    \item $c_{t,r}^j = 0$ whenever $t>r$, $r<0$, or $t<0$.
\end{itemize}

Now turn to counting paths from vertex $(n)$ in level 0 to vertex $(n-tj,tj)$ in level $r$ of $\mathcal{MN}_j(S_n,S_{n-j})$. Label a path of $r$ full-level steps ($2r$ half-level steps) starting at vertex $(n)$ in $\mathcal{MN}_j(S_n,S_{n-j})$ as an ordered $(2r+1)$-tuple of partitions, $(\lambda_0 = (n), \lambda_{\frac{1}{2}} = (n-j), \ldots, \lambda_r)$, where partitions of $n$ are indexed by integers $0,1,\ldots r$ and partitions of $n-j$ are indexed by half-integers $\half, \frac{3}{2}, \ldots \frac{2r-1}{2}$. For whole integers $i$, $\lambda_i$ is obtained from $\lambda_{i-\half}$ by the $i$th rim-$j$ hook addition, and $\lambda_{i-\half}$ is obtained from $\lambda_{i-1}$ by the $i$th rim-$j$ hook removal. Let $\mu = (n-tj,tj)$ and $P_{\mu,r}^j$ be the set of all valid paths $\{(\lambda_0,\ldots,\lambda_r): \lambda_0 = (n), \lambda_r = \mu\}$ so that $\abs{P_{\mu,r}^j} = p_{t,r}^j$.

We perform a last-step analysis. To arrive at $\lambda_r = \mu$ after the $r$th and final addition step, any path in $P_{\mu,r}^j$ must have $\lambda_{r-1}$ belonging to one of three sets:
\begin{itemize}
    \item $A_{1} = \{(n-(t-1)j, (t-1)j)\}$
    \item $A_2 = \{(n-tj,tj)\} \cup \{(n-tj,(t-1)j,j-i,1^i): 0\leq i\leq j-1\} $
    \item $A_{3} = \{(n-(t+1)j, (t+1)j)\} \cup \{(n-(t+1)j,tj,j-i,1^i): 0\leq i\leq j-1\}$.
\end{itemize}

In other words, we can express $P_{\mu,r}^j$ as a disjoint union \[P_{\mu,r}^j = \bigsqcup_{i=1}^3 P_{A_i,r}^j\] where $P_{A_i,r}^j$ is the set of paths $\{(\lambda_0,\ldots,\lambda_r):\lambda_0 = (n),\lambda_{r-1}\in A_i,\lambda_r = \mu\}$. The total path count $p_{t,r}^j$ is then equal to the sum $\abs{P_{A_1,r}^j}+\abs{P_{A_2,r}^j}+\abs{P_{A_3,r}^j}$.

Notice that $p_{t-1,r-1}^j$ counts the number of paths from $(n)$ in level 0 to the lone element of $A_1$ in level $r-1$. Furthermore, we can extend each path reaching $A_{1}$ in $r-1$ steps to connect to vertex $(n-tj,tj)$ in level $r$ via exactly 1 path by specifying $\lambda_{\frac{2r-1}{2}} = (n-tj,(t-1)j)$ and $\lambda_r = \mu$, revealing $\abs{P_{A_1,r}^j} = p_{t-1,r-1}^j$. 

We will show similarly that $\abs{P_{A_2,r}^j} = (tj+1)p_{t,r-1}^j$ and $\abs{P_{A_3,r}^j} = j(t+1)p_{t+1,r-1}^j$, yielding the desired recurrence relation. 

Notice that $A_2$ and $A_3$ have a similar structure. We prove the following claim, which is useful in analyzing the contributions to $p_{t,r}^j$ of paths through $A_2$ and $A_3$, via combinatorial argument. \\

\noindent\textbf{Claim: } For integers $1\leq t \leq r$, $n\geq 2rj$, let $\mu = (n-tj,tj)$, $P_{\mu,r}^j$ be the set of distinct paths from $(n)$ in level 0 to $\mu$ in level $r$ through $\mathcal{MN}_j(S_n,S_{n-j})$, and \[S = \{(n-tj,(t-1)j,j-i,1^i): 0\leq i\leq j-1\}.\] Denote by $P_{S,r}^j$ the set of distinct paths in $r$ full-level steps from $(n)$ in level 0 to $\lambda_r \in S$ in level $r$ through $\mathcal{MN}_j(S_n,S_{n-j})$. Then $\abs{P_{S,r}^j} = (tj-1)\abs{P_{\mu,r}^j}$. \\

\noindent\textit{Proof of Claim: } First, let us briefly describe how to identify a partition $\lambda$ with its $j$-abacus configuration, in which rim-$j$ addition and removal steps are conveniently realized as bead moves on $j$ runners.

The $j$ runners are alternately referred to as columns. Rows are indexed by integers with row labels increasing down columns. The empty partition has beads filling row 0 and all rows indexed by negative integers above. Reading the $j$-abacus diagram corresponding to integer partition $\lambda = (\lambda_1, \lambda_2, \ldots)$ from bottom to top and then right to left, $\lambda_i$ is the number of empty positions weakly above and to the left of the $i$th bead. See James and Kerber \cite[Chapter 2]{JK81} for details.

The following example serves to make the above description clear. Let $\lambda = (n-27,8,6,5,4,2^2) \vdash n$ with 3-abacus bead configuration displayed below.

\begin{figure}[H]
    \centering
    \begin{tikzpicture}[every node/.style={font=\small}]

\node at (0,0.6)   {\(\mathbf{1}\)};
\node at (1.5,0.6) {\(\mathbf{2}\)};
\node at (3,0.6)   {\(\mathbf{3}\)};

\node[anchor=east] at (-0.45,  0.0) {\(\mathbf{-1}\)};
\node[anchor=east] at (-0.45, -0.6) {\(\mathbf{0}\)};
\node[anchor=east] at (-0.45, -1.2) {\(\mathbf{1}\)};
\node[anchor=east] at (-0.45, -1.8) {\(\mathbf{2}\)};
\node[anchor=east] at (-0.45, -2.4) {\(\mathbf{3}\)};
\node[anchor=east] at (-0.45, -3.0) {\(\mathbf{4}\)};
\node[anchor=east] at (-0.45, -3.6) {\(\mathbf{5}\)};
\node[anchor=east] at (-0.45, -4.2) {\(\vdots\)};
\node[anchor=east] at (-0.45, -4.8) {\(\mathbf{\ast}\)};

\node at (0,0)      {\(\bigcirc\)};
\node at (0,-0.6)   {\(\cdot\)};
\node at (0,-1.2)   {\(\bigcirc\)};
\node at (0,-1.8)   {\(\bigcirc\)};
\node at (0,-2.4)   {\(\cdot\)};
\node at (0,-3.0)   {\(\cdot\)};
\node at (0,-3.6)   {\(\cdot\)};
\node at (0,-4.2)   {\(\vdots\)};
\node at (0,-4.8)   {\(\cdot\)};

\node at (1.5,0)      {\(\bigcirc\)};
\node at (1.5,-0.6)   {\(\cdot\)};
\node at (1.5,-1.2)   {\(\cdot\)};
\node at (1.5,-1.8)   {\(\cdot\)};
\node at (1.5,-2.4)   {\(\bigcirc\)};
\node at (1.5,-3.0)   {\(\bigcirc\)};
\node at (1.5,-3.6)   {\(\cdot\)};
\node at (1.5,-4.2)   {\(\vdots\)};
\node at (1.5,-4.8)   {\(\cdot\)};

\node at (3,0)      {\(\bigcirc\)};
\node at (3,-0.6)   {\(\bigcirc\)};
\node at (3,-1.2)   {\(\cdot\)};
\node at (3,-1.8)   {\(\bigcirc\)};
\node at (3,-2.4)   {\(\cdot\)};
\node at (3,-3.0)   {\(\cdot\)};
\node at (3,-3.6)   {\(\cdot\)};
\node at (3,-4.2)   {\(\vdots\)};
\node at (3,-4.8)   {\(\bigcirc*\)};

\end{tikzpicture}
\end{figure}

In the example above, the bead in row $\ast$ is positioned so that $\lambda_1 = n-27$ with $n$ unspecified but assumed to be at least $2rj$. In particular, the bead need not necessarily occupy the third column as in the example above. Valid rim-$j$ hook additions to $\lambda$ are the possible downward bead moves (increasing the row index of the bead by 1), and valid removal steps are performed by upward bead moves (decreasing the row index by 1). For a path $(\lambda_0 = (n), \ldots, \lambda_r = \mu)$ in $P_{\mu,r}^j$, call addition step $i$ the $k$th \textit{distinguished step} if: 
\begin{itemize}
    \item The bead with greatest row label (possibly ignoring the one labeled $\ast$) in column $j$ is moved from row $k-1$ to row $k$.
    \item That bead is never again moved to row $k-1$ in subsequent removal steps $i+1,\ldots,r$.
\end{itemize} 
Every path in $P_{\mu,r}^j$ is guaranteed to have exactly $t$ distinguished steps because the abacus corresponding to partition $(n-tj,tj)$ has a bead in row $t$ of column $j$, and no other beads with row index greater than 0 besides the one in row $\ast$. In the procedures that follow, we always ignore the bead in row $\ast$ which gives $\lambda$ its long first row of length $n-tj$, unless we say otherwise. For example, when we refer to ``the bead with greatest row index in column $g$", we do not mean the bead in row $\ast$ if it happens to be in that column.

We'll demonstrate that $\abs{P_{S,r}^j} = (tj-1)\abs{P_{\mu,r}^j}$ by showing both of the following.
\begin{enumerate}
    \item Every path in $P_{\mu,r}^j$ has exactly $tj-1$ distinct \textit{modifications} in $P_{S,r}^j$.
    \item Every path in $P_{S,r}^j$ can be realized as a modification of exactly one path in $P_{\mu,r}^j$.
\end{enumerate}

The procedure to modify a path $p$ in $P_{\mu,r}^j$ to obtain a \textit{modification} $p^\prime \in P_{S,r}^j$ is as follows.

\begin{enumerate}
    \item Select one of the $t$ distinguished addition steps in $p$. Say the $k$th distinguished addition step from $\lambda_{i_k-\half}$ to $\lambda_{i_k}$ is chosen. 
    \item Do one of the following to obtain $\lambda_{i_k}$ in $p^\prime$.
        \begin{enumerate}
            \item If $k=1$, instead of moving the bead with greatest row index in column $j$ (ignoring $\ast$), move the bead with greatest row index in one of the other $j-1$ columns, say column $g$.
            \item If $k>1$, move the bead with greatest row index in column $g \neq j$ or move the bead with the second greatest row index in column $j$ (ignoring $\ast$ again). 
        \end{enumerate}
    \item From $\lambda_{i_k}$ onward, continue making bead moves in path $p^\prime$ to match the ones made in $p$. If a bead with $d$th largest row index in column $c$ is moved in path $p$, move it also in path $p^\prime$, with the following exception. If in path $p$, the second-largest row-indexed bead in column $j$ moves downward to occupy the row directly above the largest row-indexed bead in column $j$, this move is not available in $p^\prime$. Instead move the second-largest row-indexed bead in column $g$ downward in this case. When this bead moves upward again in column $j$ in path $p$, move it upward in column $g$ in path $p^\prime$.
\end{enumerate}

See that for every path $p \in P_{\mu,r}^j$, there are exactly $tj - 1$ distinct modifications of $p$ that can be obtained by the procedure above, $j$ of them for each of the $t-1$ distinguished steps after the first one, and $j-1$ of them by modifying the first distinguished step.

This procedure is reversible in the sense that a path $p^\prime$ in $P_{S,r}^j$ can be uniquely identified as a modification of exactly one path $p$ in $P_{\mu,r}^j$ by doing the following.
\begin{enumerate}
    \item The partition $\lambda_r$ in $p^\prime$ belongs to set $S$, meaning that it has a bead in column $j$, row $t-1$, and another single bead with row index $1$, in column $g$ with $1\leq g\leq j$. First, look at the bead move history for the bead with largest row index in column $j$ and identify the $t-1$ distinguished moves for this bead. Then examine the bead move history for the bead with largest row index in column $g$ (or second-largest row index if $g=j$) and identify addition step $i$ as the final index-increasing move this bead experiences in $p^\prime$. 
    \item Change addition step $i$ from $\lambda_{i-\half}$ to $\lambda_i$ in $p$, the reverse-modification of $p^\prime$ so that this step moves the largest-index bead in column $j$ instead of the chosen bead in column $g$.
    \item Modify the remaining partitions $\lambda_{i+\half},\lambda_{i+1},\ldots,\lambda_r$ in path $p$ to match the bead moves in the corresponding steps of $p^\prime$. If a bead with $d$th largest row index in column $g\neq j$ would move to a row in column $g$ that is unavailable, move the corresponding $d$th-largest row indexed bead in column $j$ instead. 
\end{enumerate}

Since this reverse procedure results in exactly one path $p \in P_{\mu,r}^j$ for any path in $P_{S,r}^j$, this establishes the desired relationship, $\abs{P_{S,r}^j} = (tj-1)\abs{P_{\mu,r}^j}$.\hfill $\square$\\

\Yboxdim{6pt}  
\Ylinethick{.6pt} 
\begin{figure}[H]
\[
\begin{tikzpicture}[line width=.5pt,xscale=0.16,yscale=0.24]  
\path (-15,0)  node[anchor=west]  {\boldsymbol{$\ell = 0$}};
\path (10,0)  node[anchor=west] (S18-0) {\yng(18)};

\path (-15,-7.5)  node[anchor=west]  {\boldsymbol{$\ell = \half$}};
\path (25,-7.5)  node[anchor=west] (S15-1) {\yng(15)};

\path (-15,-15)  node[anchor=west]  {\boldsymbol{$\ell = 1$}};
\path (30,-15)  node[anchor=west] (S1521-2) {\yng(15,2,1)};

\path (-15,-22.5)  node[anchor=west]  {\boldsymbol{$\ell = \frac{3}{2}$}};
\path (15,-22.5)  node[anchor=west] (S1221-3) {\yng(12,2,1)};

\path (-15,-30)node[anchor=west]  {\boldsymbol{$\ell = 2$}};
\path (-3,-30)  node[anchor=west] (S1233-4) {\yng(12,3,3)};

\path (-15,-37.5)  node[anchor=west]  {\boldsymbol{$\ell = \frac{5}{2}$}};
\path (15,-37.5)  node[anchor=west] (S123-5) {\yng(12,3)};

\path (-15,-45)node[anchor=west]  {\boldsymbol{$\ell = 3$}};
\path (-3,-45)  node[anchor=west] (S126-6) {\yng(12,6)};

\path  (S18-0) edge[black,thick] (S15-1);
\path  (S1521-2) edge[red,thick] (S15-1);
\path  (S1521-2) edge[black,thick] (S1221-3);
\path  (S1233-4) edge[red,thick] (S1221-3);
\path  (S1233-4) edge[black,thick] (S123-5);
\path  (S123-5) edge[black,thick] (S126-6);

\end{tikzpicture}
\]
\caption{A path from $(18)$ in level 0 to $\mu = (12,6)$ in level 3 through $\mathcal{MN}_3(S_{18},S_{15})$, a member of $P_{\mu,3}^3$. Here $r=j=3$ and $t=2$. The corresponding $3$-abacus diagrams in Figure \ref{fig:abpath1} reveal that addition steps 2 and 3 are distinguished.\label{fig:path1}}
\end{figure}

\begin{figure}[H]
    \centering
    \begin{tikzpicture}[>=Latex, every node/.style={font=\small},xscale=0.5,yscale=0.6]

        \begin{scope}[shift={(0,0)}]
            \node at (0,0.6)   {\(\mathbf{1}\)};
            \node at (1.5,0.6) {\(\mathbf{2}\)};
            \node at (3,0.6)   {\(\mathbf{3}\)};

            \node[anchor=east] at (-0.45,  0.0) {\(\mathbf{-1}\)};
            \node[anchor=east] at (-0.45, -0.6) {\(\mathbf{0}\)};
            \node[anchor=east] at (-0.45, -1.2) {\(\mathbf{1}\)};
            \node[anchor=east] at (-0.45, -1.8) {\(\mathbf{2}\)};
            \node[anchor=east] at (-0.45, -2.4) {\(\mathbf{3}\)};
            \node[anchor=east] at (-0.45, -3.0) {\(\mathbf{4}\)};
            \node[anchor=east] at (-0.45, -3.6) {\(\mathbf{5}\)};
            \node[anchor=east] at (-0.45, -4.2) {\(\mathbf{6}\)};
            \node[anchor=east] at (-0.45, -4.8) {\(\mathbf{7}\)};

            \node at (0,0)      {\(\bigcirc\)};
            \node at (0,-0.6)   {\(\bigcirc\)};
            \node at (0,-1.2)   {\(\cdot\)};
            \node at (0,-1.8)   {\(\cdot\)};
            \node at (0,-2.4)   {\(\cdot\)};
            \node at (0,-3.0)   {\(\cdot\)};
            \node at (0,-3.6)   {\(\cdot\)};
            \node at (0,-4.2)   {\(\cdot\)};
            \node at (0,-4.8)   {\(\bigcirc\)};

            \node at (1.5,0)      {\(\bigcirc\)};
            \node at (1.5,-0.6)   {\(\bigcirc\)};
            \node at (1.5,-1.2)   {\(\cdot\)};
            \node at (1.5,-1.8)   {\(\cdot\)};
            \node at (1.5,-2.4)   {\(\cdot\)};
            \node at (1.5,-3.0)   {\(\cdot\)};
            \node at (1.5,-3.6)   {\(\cdot\)};
            \node at (1.5,-4.2)   {\(\cdot\)};
            \node at (1.5,-4.8)   {\(\cdot\)};

            \node at (3,0)      {\(\bigcirc\)};
            \node at (3,-0.6)   {\(\bigcirc\)};
            \node at (3,-1.2)   {\(\cdot\)};
            \node at (3,-1.8)   {\(\cdot\)};
            \node at (3,-2.4)   {\(\cdot\)};
            \node at (3,-3.0)   {\(\cdot\)};
            \node at (3,-3.6)   {\(\cdot\)};
            \node at (3,-4.2)   {\(\cdot\)};
            \node at (3,-4.8)   {\(\cdot\)};
        \end{scope}

        \draw[->, thick, shorten >=10pt, shorten <=10pt]
            (3.25,-2.4) -- node[midway, above] {$\frac{1}{2}$} (8.05,-2.4);

        \begin{scope}[shift={(8.6,0)}]
            \node at (0,0.6)   {\(\mathbf{1}\)};
            \node at (1.5,0.6) {\(\mathbf{2}\)};
            \node at (3,0.6)   {\(\mathbf{3}\)};

            \node at (0,0)      {\(\bigcirc\)};
            \node at (0,-0.6)   {\(\bigcirc\)};
            \node at (0,-1.2)   {\(\cdot\)};
            \node at (0,-1.8)   {\(\cdot\)};
            \node at (0,-2.4)   {\(\cdot\)};
            \node at (0,-3.0)   {\(\cdot\)};
            \node at (0,-3.6)   {\(\cdot\)};
            \node at (0,-4.2)   {\(\color{blue}\bigcirc\)};
            \node at (0,-4.8)   {\(\cdot\)};

            \node at (1.5,0)      {\(\bigcirc\)};
            \node at (1.5,-0.6)   {\(\bigcirc\)};
            \node at (1.5,-1.2)   {\(\cdot\)};
            \node at (1.5,-1.8)   {\(\cdot\)};
            \node at (1.5,-2.4)   {\(\cdot\)};
            \node at (1.5,-3.0)   {\(\cdot\)};
            \node at (1.5,-3.6)   {\(\cdot\)};
            \node at (1.5,-4.2)   {\(\cdot\)};
            \node at (1.5,-4.8)   {\(\cdot\)};

            \node at (3,0)      {\(\bigcirc\)};
            \node at (3,-0.6)   {\(\bigcirc\)};
            \node at (3,-1.2)   {\(\cdot\)};
            \node at (3,-1.8)   {\(\cdot\)};
            \node at (3,-2.4)   {\(\cdot\)};
            \node at (3,-3.0)   {\(\cdot\)};
            \node at (3,-3.6)   {\(\cdot\)};
            \node at (3,-4.2)   {\(\cdot\)};
            \node at (3,-4.8)   {\(\cdot\)};
        \end{scope}

        \draw[->, thick, shorten >=6pt, shorten <=6pt]
            (10.1,-4.85) -- node[midway, left] {$1$} (10.1,-7.45);

        \begin{scope}[shift={(8.6,-8.2)}]

            \node at (0,0)      {\(\bigcirc\)};
            \node at (0,-0.6)   {\(\bigcirc\)};
            \node at (0,-1.2)   {\(\cdot\)};
            \node at (0,-1.8)   {\(\cdot\)};
            \node at (0,-2.4)   {\(\cdot\)};
            \node at (0,-3.0)   {\(\cdot\)};
            \node at (0,-3.6)   {\(\cdot\)};
            \node at (0,-4.2)   {\(\bigcirc\)};

            \node at (1.5,0)      {\(\bigcirc\)};
            \node at (1.5,-0.6)   {\(\cdot\)};
            \node at (1.5,-1.2)   {\(\color{blue}\bigcirc\)};
            \node at (1.5,-1.8)   {\(\cdot\)};
            \node at (1.5,-2.4)   {\(\cdot\)};
            \node at (1.5,-3.0)   {\(\cdot\)};
            \node at (1.5,-3.6)   {\(\cdot\)};
            \node at (1.5,-4.2)   {\(\cdot\)};

            \node at (3,0)      {\(\bigcirc\)};
            \node at (3,-0.6)   {\(\bigcirc\)};
            \node at (3,-1.2)   {\(\cdot\)};
            \node at (3,-1.8)   {\(\cdot\)};
            \node at (3,-2.4)   {\(\cdot\)};
            \node at (3,-3.0)   {\(\cdot\)};
            \node at (3,-3.6)   {\(\cdot\)};
            \node at (3,-4.2)   {\(\cdot\)};
        \end{scope}

        \draw[->, thick, shorten >=10pt, shorten <=10pt]
            (8.05,-10.6) -- node[midway, above] {$\frac{3}{2}$} (3.25,-10.6);

        \begin{scope}[shift={(0,-8.2)}]

            \node[anchor=east] at (-0.45,  0.0) {\(\mathbf{-1}\)};
            \node[anchor=east] at (-0.45, -0.6) {\(\mathbf{0}\)};
            \node[anchor=east] at (-0.45, -1.2) {\(\mathbf{1}\)};
            \node[anchor=east] at (-0.45, -1.8) {\(\mathbf{2}\)};
            \node[anchor=east] at (-0.45, -2.4) {\(\mathbf{3}\)};
            \node[anchor=east] at (-0.45, -3.0) {\(\mathbf{4}\)};
            \node[anchor=east] at (-0.45, -3.6) {\(\mathbf{5}\)};
            \node[anchor=east] at (-0.45, -4.2) {\(\mathbf{6}\)};

            \node at (0,0)      {\(\bigcirc\)};
            \node at (0,-0.6)   {\(\bigcirc\)};
            \node at (0,-1.2)   {\(\cdot\)};
            \node at (0,-1.8)   {\(\cdot\)};
            \node at (0,-2.4)   {\(\cdot\)};
            \node at (0,-3.0)   {\(\cdot\)};
            \node at (0,-3.6)   {\(\color{blue}\bigcirc\)};
            \node at (0,-4.2)   {\(\cdot\)};

            \node at (1.5,0)      {\(\bigcirc\)};
            \node at (1.5,-0.6)   {\(\cdot\)};
            \node at (1.5,-1.2)   {\(\bigcirc\)};
            \node at (1.5,-1.8)   {\(\cdot\)};
            \node at (1.5,-2.4)   {\(\cdot\)};
            \node at (1.5,-3.0)   {\(\cdot\)};
            \node at (1.5,-3.6)   {\(\cdot\)};
            \node at (1.5,-4.2)   {\(\cdot\)};

            \node at (3,0)      {\(\bigcirc\)};
            \node at (3,-0.6)   {\(\bigcirc\)};
            \node at (3,-1.2)   {\(\cdot\)};
            \node at (3,-1.8)   {\(\cdot\)};
            \node at (3,-2.4)   {\(\cdot\)};
            \node at (3,-3.0)   {\(\cdot\)};
            \node at (3,-3.6)   {\(\cdot\)};
            \node at (3,-4.2)   {\(\cdot\)};
        \end{scope}

        \draw[->, thick, shorten >=6pt, shorten <=6pt]
            (1.5,-13.05) -- node[midway, right] {$2$} (1.5,-15.65);

        \begin{scope}[shift={(0,-16.4)}]

            \node[anchor=east] at (-0.45,  0.0) {\(\mathbf{-1}\)};
            \node[anchor=east] at (-0.45, -0.6) {\(\mathbf{0}\)};
            \node[anchor=east] at (-0.45, -1.2) {\(\mathbf{1}\)};
            \node[anchor=east] at (-0.45, -1.8) {\(\mathbf{2}\)};
            \node[anchor=east] at (-0.45, -2.4) {\(\mathbf{3}\)};
            \node[anchor=east] at (-0.45, -3.0) {\(\mathbf{4}\)};
            \node[anchor=east] at (-0.45, -3.6) {\(\mathbf{5}\)};

            \node at (0,0)      {\(\bigcirc\)};
            \node at (0,-0.6)   {\(\bigcirc\)};
            \node at (0,-1.2)   {\(\cdot\)};
            \node at (0,-1.8)   {\(\cdot\)};
            \node at (0,-2.4)   {\(\cdot\)};
            \node at (0,-3.0)   {\(\cdot\)};
            \node at (0,-3.6)   {\(\bigcirc\)};

            \node at (1.5,0)      {\(\bigcirc\)};
            \node at (1.5,-0.6)   {\(\cdot\)};
            \node at (1.5,-1.2)   {\(\bigcirc\)};
            \node at (1.5,-1.8)   {\(\cdot\)};
            \node at (1.5,-2.4)   {\(\cdot\)};
            \node at (1.5,-3.0)   {\(\cdot\)};
            \node at (1.5,-3.6)   {\(\cdot\)};
       
            \node at (3,0)      {\(\bigcirc\)};
            \node at (3,-0.6)   {\(\cdot\)};
            \node at (3,-1.2)   {\(\color{blue}\bigcirc\)};
            \node at (3,-1.8)   {\(\cdot\)};
            \node at (3,-2.4)   {\(\cdot\)};
            \node at (3,-3.0)   {\(\cdot\)};
            \node at (3,-3.6)   {\(\cdot\)};
        \end{scope}

        \draw[->, thick, shorten >=10pt, shorten <=10pt]
            (3.25,-18.8) -- node[midway, above] {$\frac{5}{2}$} (8.05,-18.8);

        \begin{scope}[shift={(8.6,-16.4)}]
            \node at (0,0)      {\(\bigcirc\)};
            \node at (0,-0.6)   {\(\bigcirc\)};
            \node at (0,-1.2)   {\(\cdot\)};
            \node at (0,-1.8)   {\(\cdot\)};
            \node at (0,-2.4)   {\(\cdot\)};
            \node at (0,-3.0)   {\(\cdot\)};
            \node at (0,-3.6)   {\(\bigcirc\)};

            \node at (1.5,0)      {\(\bigcirc\)};
            \node at (1.5,-0.6)   {\(\color{blue}\bigcirc\)};
            \node at (1.5,-1.2)   {\(\cdot\)};
            \node at (1.5,-1.8)   {\(\cdot\)};
            \node at (1.5,-2.4)   {\(\cdot\)};
            \node at (1.5,-3.0)   {\(\cdot\)};
            \node at (1.5,-3.6)   {\(\cdot\)};

            \node at (3,0)      {\(\bigcirc\)};
            \node at (3,-0.6)   {\(\cdot\)};
            \node at (3,-1.2)   {\(\bigcirc\)};
            \node at (3,-1.8)   {\(\cdot\)};
            \node at (3,-2.4)   {\(\cdot\)};
            \node at (3,-3.0)   {\(\cdot\)};
            \node at (3,-3.6)   {\(\cdot\)};
        \end{scope}

        \draw[->, thick, shorten >=6pt, shorten <=6pt]
            (10.1,-20.95) -- node[midway, left] {$3$} (10.1,-23.55);

        \begin{scope}[shift={(8.6,-24.6)}]

            \node[anchor=east] at (-0.45,  0.0) {\(\mathbf{-1}\)};
            \node[anchor=east] at (-0.45, -0.6) {\(\mathbf{0}\)};
            \node[anchor=east] at (-0.45, -1.2) {\(\mathbf{1}\)};
            \node[anchor=east] at (-0.45, -1.8) {\(\mathbf{2}\)};
            \node[anchor=east] at (-0.45, -2.4) {\(\mathbf{3}\)};
            \node[anchor=east] at (-0.45, -3.0) {\(\mathbf{4}\)};
            \node[anchor=east] at (-0.45, -3.6) {\(\mathbf{5}\)};

            \node at (0,0)      {\(\bigcirc\)};
            \node at (0,-0.6)   {\(\bigcirc\)};
            \node at (0,-1.2)   {\(\cdot\)};
            \node at (0,-1.8)   {\(\cdot\)};
            \node at (0,-2.4)   {\(\cdot\)};
            \node at (0,-3.0)   {\(\cdot\)};
            \node at (0,-3.6)   {\(\bigcirc\)};

            \node at (1.5,0)      {\(\bigcirc\)};
            \node at (1.5,-0.6)   {\(\bigcirc\)};
            \node at (1.5,-1.2)   {\(\cdot\)};
            \node at (1.5,-1.8)   {\(\cdot\)};
            \node at (1.5,-2.4)   {\(\cdot\)};
            \node at (1.5,-3.0)   {\(\cdot\)};
            \node at (1.5,-3.6)   {\(\cdot\)};

            \node at (3,0)      {\(\bigcirc\)};
            \node at (3,-0.6)   {\(\cdot\)};
            \node at (3,-1.2)   {\(\cdot\)};
            \node at (3,-1.8)   {\(\color{blue}\bigcirc\)};
            \node at (3,-2.4)   {\(\cdot\)};
            \node at (3,-3.0)   {\(\cdot\)};
            \node at (3,-3.6)   {\(\cdot\)};
        \end{scope}

    \end{tikzpicture}
    \caption{Sequential bead moves on a $3$-abacus corresponding to the path in Figure \ref{fig:path1}. Each abacus corresponds to a Young diagram-labeled vertex, and arrows correspond to edges terminating at the level indicated by the half-integer labels in $\mathcal{MN}_j(S_{18},S_{15})$. Beads are colored blue to emphasize they were moved in the prior addition/removal step. Writing out the 5 valid modifications of this path in $P_{S,3}^3$ is a good exercise. \label{fig:abpath1}}
\end{figure}

Returning to the problem at hand, we next want to show that $\abs{P_{A_2,r}^j} = (tj+1)p_{t,r-1}^j$. Notice that for every path from $(n)$ to $(n-tj,tj)$ in $r-1$ steps, there are two paths from $(n-tj,tj)$ to $(n-tj,tj)$ in the $r$th and final step. These two paths are determined by the choice of $\lambda_{\frac{2r-1}{2}}$ as either $(n-tj,(t-1)j)$ or $(n-(t+1)j,tj)\}$. Meanwhile there is one path from each of the other partitions in $A_2$, which have the form $(n-tj,(t-1)j,j-i,1^i)$ with $0\leq i \leq j-1$, at level $r-1$ to $(n-tj,tj)$ in level $r$ (here $\lambda_{\frac{2r-1}{2}}$ must be $(n-tj,(t-1)j)$). By the previous Claim, the total number of paths in $r-1$ steps to partitions in $\{(n-tj,(t-1)j,j-i,1^i): 0\leq i \leq j-1\}$ are $tj-1$ times as numerous as the paths to $(n-tj,tj)$ in $r-1$ steps, provided $t < r$ (noting that if $t=r$, there are no paths to $(n-tj,tj)$ in $r-1$ steps). In sum we have $tj-1+2 = tj+1$ paths from $(n)$ to $(n-tj,tj)$ through $A_2$ in $r$ steps for each path from $(n)$ to $(n-tj,tj)$ in $r-1$ steps ($p_{t,r-1}^j$) as desired.

Finally, we show $\abs{P_{A_3,r}^j} = j(t+1)p_{t+1,r-1}^j$. There is exactly one choice of path from each member of $A_3$ in level $r-1$ to $(n-tj,tj)$ in level $r$ (via $\lambda_{\frac{2r-1}{2}}=(n-(t+1)j,tj)$). Applying the Claim with integer $t+1$ in place of $t$, there are $(t+1)j-1$ paths in total to partitions in $\{(n-(t+1)j,tj,j-i,1^i): 0 \leq i \leq j-1\}$ in $r-1$ steps (provided $t+1\leq r-1$) for each one to $(n-(t+1)j,(t+1)j)$ in $r-1$ steps. Since $(n-(t+1)j,(t+1)j)$ is itself an element of $A_3$, this means there are $j(t+1)$ total paths from $(n)$ in level 0 to $(n-tj,tj)$ in level $r$ through $A_3$ for each of the $p_{t+1,r-1}^j$ paths arriving at $(n-(t+1)j,(t+1)j)$ in level $r-1$. Hence $\abs{P_{A_3,r}^j}=j(t+1)p_{t+1,r-1}^j$ (again recalling that if $t+1>r-1$, $p_{t+1,r-1}^j$ = 0). 

Thus we have the recurrence relation,
\[p_{t,r}^j = p_{t-1,r-1}^j + (tj+1)p_{t,r-1}^j + j(t+1)p_{t+1,r-1}^j,\]
with $p_{t,r}^j$ subject to the same boundary conditions as $c_{t,r}^j$, and the two sequences are equivalent.
\end{proof}

With the discussion of \textit{distinguished steps} and \textit{modifications} from the previous proof in mind, we are ready to complete the proof of Theorem \ref{thm:A}.

\begin{proof}[Proof of Theorem \ref{thm:A}]
Let \(r,n,j\) satisfy the assumptions of Theorem \ref{thm:A}. It is enough to
prove \eqref{eqn:mform} for \(0\leq t\leq r\), since every partition reachable
from \((n)\) after \(r\) full steps in
\(\mathcal{MN}_j(S_n,S_{n-j})\) has the form
\[
\lambda=(n-tj,\bar{\lambda}),
\qquad
\bar{\lambda}\vdash tj,
\]
for some \(0\leq t\leq r\).

Fix such a \(t\), set
\[
\mu=(n-tj,tj),
\]
and let \(P_{\mu,r}^j\) denote the set of paths in
\(\mathcal{MN}_j(S_n,S_{n-j})\) from \((n)\) at level \(0\) to
\(\mu\) at level \(r\). By Theorem \ref{thm:paths},
\[
\abs{P_{\mu,r}^j}=c_{t,r}^j.
\]
We will construct a bijection
\[
P_{\mu,r}^j\times\mathcal{T}_j(\bar{\lambda})
\longrightarrow
P_{\lambda,r}^j,
\]
where \(P_{\lambda,r}^j\) is the set of paths from \((n)\) to
\(\lambda=(n-tj,\bar{\lambda})\), and
\(\mathcal{T}_j(\bar{\lambda})\) is the set of rim-\(j\)-hook tableaux of
shape \(\bar{\lambda}\) and content \((j^t)\).

Represent each path as a sequence of \(2r\) bead moves on a \(j\)-abacus. On
each runner, label the beads according to their order from bottom to top and
preserve these labels as the beads move. Thus the history of each individual
bead is well defined. As before, we ignore the bead determining the long first
row unless explicitly stated otherwise.

Let
\[
p=(\lambda_0=(n),\lambda_{\frac12},\lambda_1,\ldots,
\lambda_{r-\frac12},\lambda_r=\mu)
\in P_{\mu,r}^j.
\]
The path \(p\) has exactly \(t\) distinguished addition steps
\[
i_1<\cdots<i_t.
\]
At the \(k\)th distinguished step, the bead with greatest row index in runner
\(j\) moves downward from row \(k-1\) to row \(k\), and that displacement is
not subsequently reversed.

Fix a rim-\(j\)-hook tableau
\[
T\colon
\varnothing=\bar{\lambda}^{(0)}
\subset\bar{\lambda}^{(1)}
\subset\cdots\subset
\bar{\lambda}^{(t)}=\bar{\lambda},
\]
where each
\(\bar{\lambda}^{(k)}/\bar{\lambda}^{(k-1)}\) is a rim-\(j\) hook. We modify
\(p\) to obtain a path
\[
p'=(\lambda_0'=(n),\lambda_{\frac12}',\lambda_1',\ldots,
\lambda_{r-\frac12}',\lambda_r')
\]
as follows.

At the distinguished step \(i_k\), replace the downward bead move in runner
\(j\) from row \(k-1\) to row \(k\) by the downward bead move corresponding to
the addition
\[
\bar{\lambda}^{(k-1)}
\subset\bar{\lambda}^{(k)}.
\]
At every non-distinguished step, perform the move of the corresponding labeled
bead whenever that move remains available. If, between \(i_k\) and
\(i_{k+1}\), a move of the bead with second greatest row index in runner \(j\)
is unavailable because of the modification at \(i_k\), replace it by the
corresponding move of the bead with the next greatest row index in the runner
used at step \(i_k\). Thereafter, whenever the original path moves the former
bead, the modified path moves the latter bead, until the corresponding inverse
move occurs. At that point, restore the original correspondence between the
labeled beads. The same convention applies after \(i_t\).

This procedure preserves the validity of every bead move. The matched
downward--upward moves have no net effect on the terminal abacus, while the
\(t\) distinguished moves have been replaced by the \(t\) rim-\(j\)-hook
additions prescribed by \(T\). Consequently,
\[
\lambda_r'=(n-tj,\bar{\lambda}),
\]
and hence \(p'\in P_{\lambda,r}^j\).

We next describe the inverse construction. Let
\[
p'=(\lambda_0'=(n),\lambda_{\frac12}',\lambda_1',\ldots,
\lambda_{r-\frac12}',\lambda_r'=\lambda)
\in P_{\lambda,r}^j.
\]
For each labeled bead, match its downward and upward moves in last-in-first-out
order: a downward move from row \(a-1\) to row \(a\) is paired with the first
subsequent, as-yet-unmatched upward move from row \(a\) to row \(a-1\). Call a
downward move \emph{persistent} if it remains unmatched. Equivalently, a
persistent move represents a unit downward displacement of a bead that remains
present in the terminal abacus.

Apart from the bead determining the first row, the terminal abacus differs from
the initial abacus by a total downward displacement of \(t\) rows. Therefore
\(p'\) has exactly \(t\) persistent downward moves. Let
\[
i_1<\cdots<i_t
\]
be the addition steps at which they occur. Reading these moves chronologically
determines a sequence
\[
\varnothing=\bar{\lambda}^{(0)}
\subset\bar{\lambda}^{(1)}
\subset\cdots\subset
\bar{\lambda}^{(t)}=\bar{\lambda},
\]
in which every
\(\bar{\lambda}^{(k)}/\bar{\lambda}^{(k-1)}\) is a rim-\(j\) hook. Thus the
persistent moves determine a unique tableau
\(T\in\mathcal{T}_j(\bar{\lambda})\).

We recover \(p\in P_{\mu,r}^j\) by reading the bead moves of \(p'\) in
chronological order. At the \(k\)th persistent addition step \(i_k\), replace
the corresponding move of \(p'\) by the downward move in runner \(j\) from row
\(k-1\) to row \(k\). This becomes the \(k\)th distinguished addition step of
\(p\). At every nonpersistent step, move the corresponding labeled bead
whenever the required position is available. If that position is unavailable
because of an earlier persistent move, perform instead the corresponding move
of the bead in runner \(j\) having the same relative position among the beads
not involved in the first \(k\) persistent moves. Maintain this new
correspondence until the inverse move occurs, and then restore the original
correspondence.

After all \(2r\) half-steps have been processed, every matched
downward--upward pair has returned to its initial position, while the \(t\)
persistent moves have been replaced by the successive moves
\[
0\longrightarrow1\longrightarrow\cdots\longrightarrow t
\]
of the bead with greatest row index in runner \(j\). The resulting terminal
abacus is therefore the abacus of
\[
\mu=(n-tj,tj),
\]
so the recovered path belongs to \(P_{\mu,r}^j\).

The forward and reverse procedures are mutually inverse: the persistent moves
of \(p'\) recover the tableau \(T\), and replacing them by the distinguished
moves in runner \(j\) recovers \(p\). Hence
\[
P_{\lambda,r}^j
\longleftrightarrow
P_{\mu,r}^j\times\mathcal{T}_j(\bar{\lambda})
\]
is a bijection.

By the Murnaghan--Nakayama rule, all tableaux in
\(\mathcal{T}_j(\bar{\lambda})\) have the same sign, and
\[
\chi^{\bar{\lambda}}_{(j^t)}
=
\sum_{T\in\mathcal{T}_j(\bar{\lambda})}
(-1)^{\operatorname{ht}(T)}.
\]
Moreover, the sign of the path corresponding to \((p,T)\) is
\((-1)^{\operatorname{ht}(T)}\). Since
\(\abs{P_{\mu,r}^j}=c_{t,r}^j\), Theorem \ref{thm:pathcount} therefore gives
\[
\begin{aligned}
m_{\lambda,r}^j
&=
\sum_{p\in P_{\mu,r}^j}
\sum_{T\in\mathcal T_j(\bar{\lambda})}
(-1)^{\operatorname{ht}(T)}\\
&=
\abs{P_{\mu,r}^j}
\sum_{T\in\mathcal T_j(\bar{\lambda})}
(-1)^{\operatorname{ht}(T)}\\
&=
c_{t,r}^j\chi^{\bar{\lambda}}_{(j^t)}.
\end{aligned}
\]
This proves \eqref{eqn:mform}. Finally, since
\[
a_j=\frac{1}{j}\psi_j,
\]
the asserted irreducible character decomposition of \((a_j)^r\) follows from
that of \((\psi_j)^r\) by multiplication by \(j^{-r}\).
\end{proof}

\begin{remark}
It would be nice to have a bijective proof of Theorem \ref{thm:paths}, analogous to the bijections given in Section 3.3 of \cite{BH17}.
\end{remark}

\section*{Acknowledgements}
I would like to thank Persi Diaconis for suggesting that I examine cycle type statistics beyond fixed points; Sami Assaf and Greta Panova for helpful conversations related to the work of Section \ref{subsec:mnj}; and Jason Fulman for encouraging me to pursue this work and for many helpful suggestions and conversations throughout. I would also like to thank Lucas Teyssier and an anonymous referee for each pointing out an erroneous claim on $j$-cycle mixing in the random $i$-cycle shuffle setting which appeared in a previous version of this work, and for several helpful comments and suggestions.

\bibliographystyle{habbrv} 
\bibliography{main}
\end{document}